\documentclass[11pt]{amsart}
\usepackage[utf8]{inputenc}
\usepackage[T1]{fontenc}
\usepackage[english]{babel}
\usepackage{amsmath,amssymb,amsthm,mathtools,mathrsfs}
\usepackage{graphicx,xcolor,caption}
\usepackage[shortlabels]{enumitem}
\usepackage[margin=2.5cm]{geometry}
\usepackage{tikz}
\tikzset{% styles for the pictures of graphs in Section \ref{sec:graphs}
	vA/.style={circle,draw,fill=blue!25,minimum size=6mm,inner sep=0pt,font=\small},
	vB/.style={circle,draw,fill=orange!35,minimum size=6mm,inner sep=0pt,font=\small},
	arista/.style={thick},
	noarista/.style={thick,dashed,red!70!black},
}
\usepackage[hidelinks]{hyperref}

\def\calD{\mathcal D}
\def\calE{\mathcal E}
\def\calT{\mathcal T}
\def\calW{\mathcal W}
\def\calX{\mathcal X}
\def\calY{\mathcal Y}
\def\B{\mathbb B}
\def\R{\mathbb R}
\def\C{\mathbb C}

\def\tr{\text{tr}}
\def\im{\text{im}}
\def\ran{\text{ran}}
\DeclareMathOperator{\Diag}{Diag}

\theoremstyle{plain}
\newtheorem{theorem}{Theorem}[section]
\newtheorem{proposition}[theorem]{Proposition}
\newtheorem{lemma}[theorem]{Lemma}
\newtheorem{corollary}[theorem]{Corollary}
\theoremstyle{definition}
\newtheorem{definition}[theorem]{Definition}
\newtheorem{example}[theorem]{Example}
\newtheorem{notation}[theorem]{Notation}
\newtheorem{remark}[theorem]{Remark}
\numberwithin{equation}{section}

\title[Minimal Hermitian completions]{Characterizations of minimal Hermitian completions via moments and subdifferentials}

\author[T. Bottazzi]{Tamara Bottazzi$^{*,\dagger}$}
\author[A. Varela]{Alejandro Varela$^{\ddagger,\S}$}
\thanks{$^{*}$Universidad Nacional de R\'io Negro, CITECCA, Sede Andina, (8400) S.C. de Bariloche, Argentina.}
\thanks{$^{\dagger}$Consejo Nacional de Investigaciones Cient\'ificas y T\'ecnicas (CONICET), Argentina.}
\thanks{$^{\ddagger}$Instituto Argentino de Matem\'atica ``Alberto P. Calder\'on'' (IAM-CONICET), Argentina.}
\thanks{$^{\S}$Instituto de Ciencias, Universidad Nacional de General Sarmiento (UNGS), Argentina.}
\thanks{E-mail addresses: \texttt{tbottazzi@unrn.edu.ar} (T. Bottazzi), \texttt{avarela@campus.ungs.edu.ar} (A. Varela).}


\date{}

\begin{document}

\begin{abstract}
Given a Hermitian matrix $A_0$ and a real subspace $\calW$ of Hermitian matrices, we study the matrices $A_0+W_0$, $W_0\in\calW$, of minimal spectral norm, which we call $\calW$-minimal. The main motivation is the Hermitian matrix completion problem, where $\calW$ is spanned by the Hermitian matrix units corresponding to the unspecified entries. Our main tool is the moment of a subspace $S\subseteq\C^n$ relative to $\calW$, a compact convex subset of $\R^r$, $r=\dim\calW$, which we relate to the joint numerical range and to the $k$-dimensional field of a family of matrices. We compute the subdifferential of $x\mapsto\|A_0+W(x)\|$ in terms of the moments of the eigenspaces of the extreme eigenvalues of $A_0+W(x)$, and we obtain several equivalent characterizations of $\calW$-minimality: in terms of eigenvectors, of moments, of subdifferentials and of norm one functionals that annihilate $\calW$. When the identity belongs to $\calW$, every $\calW$-minimal matrix has balanced spectrum, and minimality is equivalent to the moments of the eigenspaces of the maximum and minimum eigenvalues having nonempty intersection; we show that this characterization fails if $I\notin\calW$. We illustrate the results with several examples, and we prove that the adjacency matrix of a connected graph is minimal for its completion problem if and only if the graph is complete bipartite.
\end{abstract}

\maketitle

\section{Introduction}\label{sec:intro}

Let $M_n^h(\C)$ be the real space of $n\times n$ Hermitian matrices, endowed with the spectral norm $\|\cdot\|$, and let $\calW\subseteq M_n^h(\C)$ be a real subspace. Given $A_0\in M_n^h(\C)$, we consider the best approximation problem
\[
\min_{W\in\calW}\|A_0+W\|,
\]
and we say that $A_0+W_0$ is \emph{$\calW$-minimal} if $W_0$ attains this minimum. Our main motivation is the Hermitian matrix completion problem: if some entries of a Hermitian matrix are specified and the remaining ones are free, then, setting the free entries of $A_0$ equal to zero, finding a completion of minimal spectral norm is exactly this problem for the subspace $\calW$ spanned by the Hermitian matrix units supported on the free positions (see Section \ref{sec:setting}). Since $\|\operatorname{Re}(X)\|\leq\|X\|$, the distance from $A_0$ to $\calW$ coincides with the distance to its complex span, so this is also the minimal norm among arbitrary (not necessarily Hermitian) completions.

When $\calW$ is the subspace of real diagonal matrices, minimal Hermitian matrices were characterized in \cite{alrv} and their supports were studied in \cite{soportes}, where the moment of a subspace was introduced. These ideas were extended to diagonal compact operators in \cite{BV-2013, B-21, BV-Moment}, related to the joint numerical range in \cite{KV moment, BV-Moment} and to the subdifferential of the maximum eigenvalue and semidefinite programming in \cite{BV-SDP}. In \cite{BV-2026} the case where $\calW$ is the Hermitian part of a unital C$^*$-subalgebra of $M_n(\C)$ was considered; see also \cite{rieffel, zhang-jiang-han-matrices, zhang-jiang-compactos}. Best approximation in the spectral norm by elements of a general subspace has recently received renewed attention \cite{FLT, grover26}. The classical tools for this problem are the characterization of best approximations by extremal functionals \cite{singer, LR81} and Watson's formula for the subdifferential of the spectral norm \cite{watson}.

In this paper we consider an arbitrary real subspace $\calW\subseteq M_n^h(\C)$. Neither the identity $I$ of $M_n(\C)$ nor any product structure is assumed, which is the natural setting for completion problems: the subspace spanned by the free positions of a pattern is almost never an algebra, and it contains the identity only when the whole diagonal is free. Fix an orthonormal basis $\{B^j\}_{j=1}^r$ of $\calW$. For a subspace $S\subseteq\C^n$, the moment of $S$ is the compact convex set
\[
m_S=\operatorname{conv}\left\{\left(\langle B^1v,v\rangle,\dots,\langle B^rv,v\rangle\right):\ v\in S,\ \|v\|=1\right\}\subseteq\R^r,
\]
that is, the convex hull of the field of $\{B^j\}_{j=1}^r$ restricted to $S$ in the sense of \cite{FLT}.

The absence of the identity is not a technical detail. If $I\in\calW$, every $\calW$-minimal matrix has balanced spectrum, $\lambda_{\max}=-\lambda_{\min}$, and minimality is detected by a witness that splits into two density matrices of equal weight, supported on the eigenspaces of $\lambda_{\max}$ and $\lambda_{\min}$ and with the same moment \cite{BV-2026}. For a general subspace neither property holds: the spectrum of a minimal matrix need not be balanced, and even when it is, the two parts of a witness may carry different weights. The main contribution of this paper is a description of $\calW$-minimality that accounts for both phenomena, and whose central object is the subdifferential of the norm along $\calW$, expressed through moments.

More precisely, the following results are new.
\begin{itemize}
	\item \emph{Subdifferential of the norm along $\calW$} (Theorem \ref{teo subdif norma} and Corollary \ref{formula for subdif composition norm}). For every real subspace $\calW$ and every $x\in\R^r$,
	\[
	\partial(\|A(\cdot)\|)(x)=\operatorname{conv}\left(m_{S_{\|A(x)\|}}\cup-m_{S_{-\|A(x)\|}}\right),\qquad A(x)=A_0+\sum_{j=1}^rx_jB^j,
	\]
	where $S_{\pm\|A(x)\|}$ are the eigenspaces of $A(x)$ for $\pm\|A(x)\|$ (one of them may be trivial). This single formula covers the balanced and the non-balanced cases, and reduces the computation of the subdifferential of the spectral norm, restricted to an affine subspace of Hermitian matrices, to the geometry of two compact convex subsets of $\R^r$.
	\item \emph{Characterization of $\calW$-minimality without the identity} (Theorem \ref{teo minimality merged}). $A_0+W_0$ is $\calW$-minimal if and only if $\overline 0\in m_{S_{\mathrm{MAX}}}$ when the spectrum is not balanced, and if and only if $\overline 0\in\operatorname{conv}\left(m_{S_{\lambda_{\max}}}\cup-m_{S_{\lambda_{\min}}}\right)$ when it is. We also give the equivalent formulations in terms of finitely many eigenvectors (at most $r+1$), of the subdifferential, and of functionals in the annihilator of $\calW$.
	\item \emph{The intersection criterion fails without the identity} (Remark \ref{rem:sin identidad}). When $I\in\calW$ the condition above reduces to $m_{S_{\lambda_{\max}}}\cap m_{S_{\lambda_{\min}}}\neq\emptyset$, but we give minimal matrices with balanced spectrum, for subspaces not containing $I$, for which $m_{S_{\lambda_{\max}}}\cap m_{S_{\lambda_{\min}}}=\emptyset$ and $m_{S_{\lambda_{\max}}}\cap-m_{S_{\lambda_{\min}}}=\emptyset$. One of them comes from a completion pattern (Example \ref{ex:2x2 block}).
	\item \emph{Moments and fields} (Proposition \ref{relation sets F} and Corollary \ref{cor field y momento}). We relate the moments to the $k$-dimensional field of \cite{FLT}. For $A=A_0+W_0$, the field of $\{AB^j\}_{j=1}^r$ on the subspace where $A$ attains its norm is computed in terms of the fields of $\{B^j\}_{j=1}^r$ on the eigenspaces of $\lambda_{\max}$ and $\lambda_{\min}$.
	\item \emph{Adjacency matrices of graphs} (Proposition \ref{prop:grafos}). For the completion problem in which the edges of a connected graph are specified and the non-edges and the diagonal are free, the adjacency matrix is minimal if and only if the graph is complete bipartite.
\end{itemize}

Some of the remaining results extend known ones to the present setting, and we include them, with proofs when these are short, to keep the paper self-contained. The description of the moment through density matrices and the joint numerical range (Proposition \ref{prop momento}) follows \cite{BV-2026}. The characterization by a single witness (Theorem \ref{teo AX modulo}) is \cite[Theorem 2.6]{zhang-jiang-han-matrices} and when $I\in\calW$, the equivalence between minimality and the intersection of the moments (Theorem \ref{teo mini equiv general}) was proved in \cite{BV-2026} for C$^*$-subalgebras. Moreover, the support conditions of the witness come from \cite{alrv}. We observe that these proofs use $\calW$ only through its annihilator and through $I\in\calW$, so they hold for arbitrary subspaces containing the identity. We also describe the extremal functionals on $M_n^h(\C)$ (Proposition \ref{prop:extremales}), a fact that is probably known but for which we could not find a reference.

The paper is organized as follows. In Section \ref{sec:setting} we set up the completion problem, the annihilator of $\calW$, the duality formula and the subspace $S_{\mathrm{MAX}}$ where a Hermitian matrix attains its norm. Section \ref{sec:moment} is devoted to the moment of a subspace, the joint numerical range and the field. In Section \ref{sec:subdiff} we compute the subdifferentials, and in Section \ref{sec:characterization} we prove the characterizations of $\calW$-minimality. Section \ref{sec:examples} contains the examples and the application to adjacency matrices of graphs.

\section[Hermitian matrix completion and $\calW$-minimal matrices]{Hermitian matrix completion and\\ $\calW$-minimal matrices}\label{sec:setting}
Let $A_0\in M_n^h(\C)$ and $I_0\subset \{(i,j):\ 1\leq i, j\leq n\}$ with the following property
\[(A_0)_{ij}=0,\ \text{ if }\ (i,j)\in I_0.\]
Consider the canonical basis $\{e_i\}_{i=1}^n$ of $\C^n$ and define the following $n\times n$ Hermitian matrices
\begin{equation}\label{def mat Wij}
	\begin{split}
		B^{j,j}&=e_je_j^*, \text{ for } (j,j)\in I_0,
		\\
		B^{i,j}&=\frac{1}{\sqrt{2}}\left( e_{i}e_j^* +e_{j}e_{i}^*\right), \text{ for } (i,j)\in I_0 \text{ with } i<j,	\\
		B^{i,j}&=	\frac{\mathrm{i}}{\sqrt{2}}\left( e_{i} e_j^* - e_{j} e_{i}^*\right), \text{ for }(i,j)\in I_0 \text{ with } j<i,
	\end{split}
\end{equation}
and the orthonormal set
\begin{equation}\label{base W}
	\B=\{B^{i,j}:(i,j)\in I_0\}.
\end{equation}
Consider $\calW=\operatorname{span}_{\R}\B$, which is a real subspace of $M_n^h(\C)$ with $\dim(\calW)= \# I_0=r$. Consider an order fixed in $I_0$, $\sigma: I_0\to \{1,2,\dots, r\}$ such that $B^{i,j}=B^{\sigma(i,j)}=B^k$, $1\leq k\leq r$. We consider the affine matrix
\begin{equation}\label{A(x)}
	A(x)=A_0+\sum_{k=1}^rx_{k}B^k=A_0+W(x),
\end{equation}
with $x=(x_1,x_2,\dots,x_r)\in \R^{r}$. We shall also use the linear map $\Phi^\B:M_n^h(\C)\to \R^r$ defined by
\begin{equation}\label{defi phi general}
	\Phi^\B(R)=(\tr(RB^1),\tr(RB^2),\dots,\tr(RB^r)), \text{ for every } R\in M_n^h(\C),
\end{equation}
which we will later restrict to the set $\calD_S$ of density matrices (see \eqref{defi phi}).
In this context, we state the following approximation problem: 
\begin{equation}\label{completion problem}
	\min_{x\in \R^r}\|A(x)\|,
\end{equation}
where $\|\cdot\|$ is the spectral norm. Note that \eqref{completion problem} is equivalent to finding $W_0=W(x_0)\in \calW$ such that
\[\text{dist} (A_0,\calW)=\|A_0+W_0\|\leq \|A_0+W\|,\ \text{ for every } W\in \calW.\]
In this case, we say that $A(x_0)=A_0+W(x_0)$ is a $\calW$-minimal matrix and $-W(x_0)$ is a best approximation of $A_0$ in the subspace $\calW$. Observe that for every $z\in \C^r$, $\operatorname{Re}(A(z))=A_0+\operatorname{Re}( W(z))$ and
\[\|A_0+\operatorname{Re}(W(z))\|=\|\operatorname{Re}(A(z))\|\leq \|A(z)\|=\|A_0+W(z)\|.\]
This implies that
\[\min_{x\in \R^r}\|A(x)\|\leq \min_{z\in \C^r}\|A(z)\|.\]
Since $\R^r\subset\C^r$, the reverse inequality is trivial, and therefore $\text{dist}(A_0,\calW)=\text{dist}(A_0,\operatorname{span}_{\C}\B)$. 

\begin{remark} \label{I0 simetrico}
	Let $A_0\in M_n^h(\C)$ and let $I_0$ as defined before. If $(i,j)\in I_0$, then $(j,i)\in I_0$, since $(A_0)_{ji}=\overline{(A_0)_{ij}}$ for every $i,j$.
\end{remark}

Define the norm 1 elements of the annihilator of $\calW$ in the dual space of $M_n^h(\C)$ as
\begin{equation}\label{anulador}
	\calW^{\circ}_1=\{X\in M_n^h(\C):\ \tr(XW)=0, \ \forall  W\in \calW \wedge \|X\|_1=1\},
\end{equation}  
where $\|X\|_1=\tr\left( (X^*X)^{1/2}\right) $ is the trace norm. Thus $\calW^{\circ}_1$ is the set of elements of trace norm one of the annihilator $\calW^{\circ}=\{X\in M_n^h(\C):\ \tr(XW)=0, \ \forall W\in\calW\}$. It is evident that $X\in \calW^{\circ} $ if and only if $\tr(XB^{i,j})=0$ for every $(i,j)\in I_0$. This implies in particular that $\ker(\Phi^\B)=\calW^{\circ}$ and $\calW^{\circ}_1=\{X\in\ker(\Phi^\B):\ \|X\|_1=1\}$, for the map $\Phi^\B$ defined in \eqref{defi phi general}. Also, we can deduce that $\dim_{\R}\calW^{\circ}=n^2-\dim_{\R}\calW=n^2-\#I_0=n^2-r$.

\begin{proposition}\label{A0 Y comparten ceros}
Let $A_0\in M_n^h(\C)$ and let $I_0=\{(i,j): (A_0)_{ij}=0\}$. For any $Y\in M_n^h(\C)$, the following conditions are equivalent:
\begin{enumerate}
	\item $Y\in \calW^{\circ}_1$.
	\item $\|Y\|_1=1$ and $Y_{ij}=0, \text{ for all } (i,j)\in I_0.$
\end{enumerate}
\end{proposition}

\begin{proof}
Consider the orthonormal set $\mathbb{B}=\{B^{i,j}:(i,j)\in I_0\}$ with $B^{i,j}$ as in \eqref{def mat Wij}.
By Remark \ref{I0 simetrico}, $B^{i,j}\in \mathbb{B} \Leftrightarrow B^{j,i}\in \mathbb{B}$, so it is sufficient to consider the cases $i<j$ and $i=j$.  Direct calculations show that
\begin{equation}\label{Yij}
	\tr(YB^{i,j})=\frac{2}{\sqrt{2}} \operatorname{Re}(Y_{ij})\quad \text{and }\quad	\tr(YB^{j,i})=-\frac{2}{\sqrt{2}} \operatorname{Im}(Y_{ij}),
\end{equation}
for $i<j$, and
\begin{equation}\label{Yii}
	\tr(YB^{i,i})=\tr\left(Ye_ie_i^*\right)=\langle Ye_i,e_i\rangle=Y_{ii},
\end{equation}
for $i=j$. \\
$(2)\Rightarrow(1)$ Since $Y_{ij}=0$ for all $(i,j)\in I_0$, it follows from \eqref{Yij} and \eqref{Yii} that $\tr(YB^{i,j})=0$ for every $(i,j)\in I_0$, and hence $Y\in\calW^{\circ}_1$.\\
$(1)\Rightarrow(2)$ Let $Y\in\calW^{\circ}_1$. Then
\[\tr(YB^{i,j})= 0, \text{ for every }\ (i,j)\in I_0.\]
By \eqref{Yij} and \eqref{Yii}, it follows that for every $(i,j)\in I_0$:
\[\left\lbrace \begin{array}{lll}
	 Y_{ij}=\overline{Y_{ji}}=0, &\text{for}&i<j \\
	Y_{ii}=0&\text{for}&i=j 
\end{array}\right. \]
Thus, $(i,j)\in I_0$ implies $Y_{ij}=0$ for all such indices.
\end{proof}

\begin{remark}[About the best approximations of real symmetric matrices]
	Let $X\in M_n(\C)$, then $\Re(X)= \frac{X+\overline{X}}{2}$ is a real matrix and $\|\Re(X)\|\leq \|X\|$. Suppose that $X=A_0+W_0$ is a $\calW$-minimal matrix with $(A_0)_{ij}\in \R$ for all $1\leq i,j\leq n$. Then,
	\[\|A_0+\Re(W_0)\|=\|\Re(A_0+W_0)\|\leq \|A_0+W_0\|.\]
	Therefore, for real symmetric matrices $A_0$ we can consider real symmetric best approximants in $\calW$. This reduces the size of the basis $\mathbb{B}$, since it contains only the elements $(B_{k,j})$ with ($k \leq j$), as defined in \eqref{def mat Wij}.
	In this case, the annihilator of $\calW$ in the dual space of $M_n^h(\R)$ is
	\[\calW^{\circ}_1=\{X\in M_n^{sim}(\R):\ \tr(XW)=0, \ \forall  W\in \calW \wedge \|X\|_1=1\}.
	\]
\end{remark}

We state the well-known Banach duality formula relating the distance from a matrix to a subspace with the norms of functionals in the dual space.
\begin{proposition}[Banach Duality Formula]
	Let $A_0\in M_n^h(\C)$. Then, 
	\begin{equation}\label{BDF}
		\max_{X\in \calW^{\circ}_1}|\tr(XA_0)|=\min_{x\in \R^r}\|A(x)\|.
	\end{equation}
\end{proposition}

Since $A_0$ vanishes on $I_0$, Proposition \ref{A0 Y comparten ceros} gives $\frac{A_0}{\|A_0\|_1}\in\calW_1^\circ$ whenever $A_0\neq0$, and the duality formula \eqref{BDF} immediately yields the following lower bound for the norm of a $\calW$-minimal matrix.

\begin{corollary}
	For every nonzero $A_0\in M_n^h(\C)$,
	\[
	\inf_{W\in\calW}\|A_0+W\|
	\ge
	\frac{\|A_0\|_F^2}{\|A_0\|_1},
	\]
	where $\|\cdot\|_F$ denotes the Frobenius norm, that is $\|X\|_F=\sqrt{\operatorname{tr}(X^*X)}$, for $X\in M_n(\C)$.
\end{corollary}
\begin{proof}
	By \eqref{BDF},
	\[
	\inf_{W\in\calW}\|A_0+W\|
	=
	\max_{X\in\calW_1^\circ}
	\left|\operatorname{tr}(XA_0)\right|
	\ge
	\left|
	\operatorname{tr}\!\left(\frac{A_0}{\|A_0\|_1}A_0\right)
	\right|
	=
	\frac{\operatorname{tr}(A_0^2)}{\|A_0\|_1}
	=
	\frac{\|A_0\|_F^2}{\|A_0\|_1}.
	\]
\end{proof}

We now introduce a subspace that plays a central role, and establish some of its properties.

\begin{proposition}[Properties of $S_{\mathrm{MAX}}$ for Hermitian matrices]\label{lema prop smax}
	Let $A\in M_n^h(\C)$, $A\neq0$, and for $\lambda\in\mathbb{R}$ let
	$S_\lambda=\{v\in\C^n:\ Av=\lambda v\}$ (so $S_\lambda=\{0\}$ if
	$\lambda\notin\sigma(A)$). Define the subspace
	\begin{equation}\label{SMAX generico}
		S_{\mathrm{MAX}}(A)=\operatorname{span}\{v\in \C^n:\ \exists\, u \in \C^n
		\text{ such that } Av=\|A\|u,\ \|v\|=\|u\|=1\}.
	\end{equation}
	Then,
	\begin{enumerate}
		\item $\begin{aligned}[t]
			S_{\mathrm{MAX}}(A)&=\operatorname{span}\{v\in\C^n:\ Av=\|A\|v
			\text{ or } Av=-\|A\|v\}\\
			&=S_{\|A\|}\oplus S_{-\|A\|}.
		\end{aligned}$
		\item Exactly one of the following conditions holds
		\begin{itemize}
			\item $S_{-\|A\|}=\{0\}$ and $S_{\mathrm{MAX}}(A)=S_{\|A\|}$,
			\item $S_{\|A\|}=\{0\}$ and $S_{\mathrm{MAX}}(A)=S_{-\|A\|}$,
			\item $S_{\|A\|}\neq\{0\}\neq S_{-\|A\|}$ and
			$S_{\mathrm{MAX}}(A)=S_{\|A\|}\oplus S_{-\|A\|}$.
		\end{itemize}
		\item $S_{\mathrm{MAX}}(A)^{\perp}=\bigoplus\limits_{\substack{\lambda\in\sigma(A)\\
				|\lambda|<\|A\|}}S_{\lambda}$.
		\item If $M_{A}$ denotes the set of unit vectors at which $A$ attains its
		norm, that is,
		\[
		M_A=\{v\in \C^n:\ \|v\|=1,\ \|Av\|=\|A\|\},
		\]
		then
		\begin{equation}\label{M_A=Esfera unitaria de S_MAX}
			M_A=S_{\mathrm{MAX}}(A)\cap \{x\in \C^n:\ \|x\|=1\}.
		\end{equation}
		\item Every unit vector $v\in S_{\mathrm{MAX}}(A)$, written as
		$v=v_++v_-$ with $v_\pm\in S_{\pm\|A\|}$, satisfies $Av=\|A\|u$ with
		\[
		u=v_+-v_-\in S_{\mathrm{MAX}}(A) \quad\text{and}\quad \|u\|=1.
		\]
	\end{enumerate}
\end{proposition}
\begin{proof}
	For $\lambda\in\sigma(A)$, let $P_\lambda$ be the orthogonal projection onto
	$S_\lambda$. By the spectral theorem, $\C^n=\bigoplus_{\lambda\in\sigma(A)}
	S_\lambda$ is an orthogonal decomposition and
	$\|A\|=\max\{|\lambda|:\lambda\in\sigma(A)\}$. Hence, for every $v\in\C^n$,
	\begin{equation}\label{eq norma Av}
		\|Av\|^2=\sum_{\lambda\in\sigma(A)}\lambda^2\|P_\lambda v\|^2
		\le\|A\|^2\sum_{\lambda\in\sigma(A)}\|P_\lambda v\|^2=\|A\|^2\|v\|^2,
	\end{equation}
	with equality if and only if $P_\lambda v=0$ whenever $|\lambda|<\|A\|$, that
	is, if and only if $v\in S_{\|A\|}+S_{-\|A\|}$. This sum is direct (and
	orthogonal) because $A\neq0$ implies $\|A\|\neq-\|A\|$.
	\begin{enumerate}
		\item The second equality holds because the span of the union of two
		subspaces is their sum, which is direct here, as noted above.
		If $\|v\|=1$ and $v\in S_{s\|A\|}$ with $s\in\{1,-1\}$, then
		$Av=\|A\|u$ with $u=sv$ and $\|u\|=1$, so $v$ is one of the generators in
		\eqref{SMAX generico}. Hence
		$S_{\|A\|}\oplus S_{-\|A\|}\subseteq S_{\mathrm{MAX}}(A)$.
		Conversely, if $\|v\|=\|u\|=1$ and $Av=\|A\|u$, then
		$\|Av\|=\|A\|=\|A\|\|v\|$, so equality holds in \eqref{eq norma Av} and
		$v\in S_{\|A\|}\oplus S_{-\|A\|}$. Since the latter is a subspace,
		$S_{\mathrm{MAX}}(A)\subseteq S_{\|A\|}\oplus S_{-\|A\|}$.
		\item Since $\|A\|=\max\{|\lambda|:\lambda\in\sigma(A)\}$, at least one of
		$\pm\|A\|$ is an eigenvalue of $A$, so $S_{\|A\|}$ and $S_{-\|A\|}$ are
		not both zero. Hence exactly one of the three listed cases occurs, and in
		each of them the stated description of $S_{\mathrm{MAX}}(A)$ follows from
		item~1.
		\item By item~1 and the orthogonal decomposition
		$\C^n=\bigoplus_{\lambda\in\sigma(A)}S_\lambda$, the orthogonal complement
		of $S_{\mathrm{MAX}}(A)$ is the sum of the eigenspaces $S_\lambda$ with
		$\lambda\in\sigma(A)\setminus\{\|A\|,-\|A\|\}$, that is, with
		$|\lambda|<\|A\|$.
		\item If $v\in S_{\mathrm{MAX}}(A)$ and $\|v\|=1$, then by item~1 equality
		holds in \eqref{eq norma Av}, so $\|Av\|=\|A\|$ and $v\in M_A$.
		Conversely, if $v\in M_A$, then $u=Av/\|A\|$ satisfies $\|u\|=1$ and
		$Av=\|A\|u$, so $v\in S_{\mathrm{MAX}}(A)$ by definition.
		\item Let $v\in S_{\mathrm{MAX}}(A)$ with $\|v\|=1$. By item~1,
		$v=v_++v_-$ with $v_\pm\in S_{\pm\|A\|}$, and
		$\|v_+\|^2+\|v_-\|^2=\|v\|^2=1$ because $v_+\perp v_-$. Then
		\[
		Av=\|A\|v_+-\|A\|v_-=\|A\|u,\qquad u=v_+-v_-\in S_{\mathrm{MAX}}(A),
		\]
		and $\|u\|^2=\|v_+\|^2+\|v_-\|^2=1$.
	\end{enumerate}
\end{proof}

\section{Moment and joint numerical range}\label{sec:moment}
Let 
\begin{equation}
	\label{def matrices de densidad}
	\calD=\{\rho\in M^h_n(\C): \rho\geq 0 \text{ and } \tr(\rho)=1\}
\end{equation}
be the set of density matrices. We now define two sets relevant to our analysis.
\begin{definition}
	Let $S$ be a subspace of $\C^n$ and $P_S$ be the orthogonal projection onto $S$.	Consider an order fixed in $I_0$, $\sigma: I_0\to \{1,2,\dots, r\}$ such that $B^{i,j}=B^{\sigma(i,j)}=B^k$, $1\leq k\leq r$, where $B^{i,j}$ are the matrices described in \eqref{def mat Wij}. Define the following family of matrices
	\[P_S\B P_S=\{P_SB^1 P_S, P_SB^2 P_S,\cdots,P_SB^r P_S\}.\]
	The joint (algebraic) numerical range with respect to the family $P_S\B P_S$ is
	\begin{equation}\label{def jnr}
		W(P_S\B P_S)= \{(\tr(\rho P_SB^1 P_S),\cdots, \tr(\rho P_SB^r P_S))\in \R^r:\ \rho\in \calD\}.
	\end{equation} 
	
\end{definition}

\begin{definition}\label{momento}
	Given a subspace $S$ of $\C^n$, we define the moment of $S$ relative to $\calW$ in the basis $\B$ as
	\begin{equation}\label{eq def momento}
		m_S^{\B}=m_S=\operatorname{conv}\{\left(\langle B^{1}v,v\rangle, \cdots, \langle B^{r}v,v\rangle\right): v\in S, \|v\|=1 \}.
	\end{equation}
\end{definition}
Some properties of $m_S$ are summarized below.
\begin{remark} \label{obs sobre momento}
	\begin{itemize}
		\item $m_S\subseteq \R^r$, since $B^i$ is Hermitian for every $1\leq i\leq r$. 
		\item The moment of a subspace was first defined in \cite{soportes} for $\calW$ as the subspace of diagonal real matrices; the name comes from the fact that it is the image of a certain moment map. In \cite{BV-Moment} it was generalized to compact operators. Recently, in \cite{BV-2026}, this was generalized by considering $\calW$ to be any unital $C^*$ subalgebra of $M_{n}(\mathbb{C})$, with identity $I = I_n$.
		\item If $\calW=D_n(\R)$ the subspace of $n\times n$ real diagonal matrices and $\B=\mathbb{E}=\{E_i=e_ie_i^*:1\leq i\leq n\}$ is the diagonal canonical basis, then
		\begin{align*}
			m_S^{\mathbb{E}}&=\operatorname{conv}\{\left(\langle E_1v,v\rangle, \cdots, \langle E_nv,v\rangle\right): v\in S, \|v\|=1 \}\\
			&=\operatorname{conv}\{\left(|v_1|^2, \cdots, |v_n|^2\right): v=(v_1,\dots,v_n)\in S, \|v\|=1 \}\\
			&=\operatorname{conv}\{|v|^2: v\in S, \|v\|=1 \},
		\end{align*}
		which is exactly the definition of moment in \cite{KV moment, BV-Moment}. 
		\item Recall that, for a family of matrices $A^1,\dots,A^k\in M_n(\C)$ and a subspace $S\subseteq\C^n$, the $k$-dimensional field of the family restricted to $S$ was defined in \cite[Definition 10]{FLT} as
		\begin{equation}\label{field}
			\mathscr{F}(\{A^j\}_{j=1}^k,S)=\{(\langle A^1v,v\rangle,\dots,\langle A^kv,v\rangle):\ v\in S,\ \|v\|=1\}\subseteq\C^k .
		\end{equation}
		In particular, since the matrices $B^j$ are Hermitian,
		\[
		\mathscr{F}(\{B^j\}_{j=1}^r,S)=\{(\langle B^1v,v\rangle,\dots,\langle B^rv,v\rangle):\ v\in S,\ \|v\|=1\}\subseteq\R^r ,
		\]
		which is exactly the set whose convex hull defines $m_S$ in \eqref{eq def momento}. Hence
		\[m_S=\operatorname{conv}\mathscr{F}(\{B^j\}_{j=1}^r,S).\]
		\item If we consider the joint numerical range with respect to the family $P_S\B P_S$ \eqref{def jnr}, 
		then $m_S\subseteq W(P_S\B P_S)$, since
		\[\langle B^{i}v,v\rangle= \tr(vv^*B^i)= \tr(P_Svv^*P_SB^i) \text{ and  } vv^*\in \calD\]
		for every $1\leq i\leq r$ and $v\in S$.
	\end{itemize}
	
\end{remark}
Let us now define the subset of density matrices relative to the subspace $S\subset \C^n$
\begin{equation}
	\label{def matrices de densidad de S}
	\calD_{S}=\{\rho\in \calD: \im(\rho)\subset S\}\subseteq \calD.
\end{equation}
Then every $\rho\in \mathcal{D}_S$ can be written as
\begin{equation}\label{equ densidad escrita usando autovals y autovecs}
	\rho=\sum_{i=1}^{\dim(S)} \lambda_i u_i u_i^*
\end{equation}
for $\lambda_i\in\R_{\geq 0}$, with $\sum_{i=1}^{\dim(S)} \lambda_i =1$ and $\{u_i\}_{i=1}^{\dim(S)}$ an orthonormal basis of $S$. 

We now restrict $\Phi^\B$ (defined in \eqref{defi phi general}) to $\calD_{S}$; that is, for every $\rho\in\calD_S$,
\begin{align}\label{defi phi}
	\begin{split}
		\Phi^{\B}(\rho)&= (\langle \rho,B^1\rangle,\langle \rho,B^2\rangle, \dots, \langle \rho,B^r\rangle)\\
		&= (\tr(\rho B^1),\tr(\rho B^2), \dots, \tr( \rho B^r)),
	\end{split}
\end{align}
where the $k$-th coordinate is the coefficient of $B^k$ in the orthogonal projection of $\rho$ onto $\calW$. Additionally, 
\[\ker(\Phi^{\B}|_{\calD_S})=\{\rho \in \calD_{S}:\ \rho\in \calW^{\circ}_1\}=\calD_S\cap \calW^{\circ}_1.\]

The following result is based on ideas from \cite[Proposition 1]{BV-2026}. For completeness, we include the proofs.

\begin{proposition} \label{prop momento}
	Let $S\subseteq \C^n$ be a subspace with $\dim(S)>0$, $\B=\{B^k\}_{k=1}^{r}$ be the fixed orthonormal basis of Hermitian matrices of a subspace $\calW\subset M_n^h(\C)$
	and $\Phi^{\mathbb{B}}:\calD_{S}\to\R^r$ defined as in \eqref{defi phi}. Then the moment $m_S$ defined in \eqref{eq def momento} admits the following descriptions.
	
	\begin{enumerate}
		\item \label{mS = diag DsubS} $m_S= \Phi^{\mathbb{B}}(\calD_{S})$.
		\item \label{item3Prop1}
		$\displaystyle{
			m_S=\{(\tr(\mu B^1),\dots,\tr( \mu B^r )): \mu\in \mathcal{D}_S\} \ \ \text{(see \eqref{def matrices de densidad de S})}
		}$
		\item \label{item1Prop1}
		$\displaystyle{
			m_S=\mathop{\bigcup}\limits_{\substack{\text{o.n. basis}  \\
					\{z_k\}_{k=1}^{\dim(S)} \text{of } S}}
			{\rm conv} \left\{\left(\langle B^1 z_k, z_k\rangle ,
			\dots,
			\langle B^r   z_k, z_k\rangle \right) \right\}_{k=1}^{\dim(S)}
		}$.
		\item \label{item4Prop1}
		If $S\neq \C^n$, then
		$W(P_SB^1P_S,\dots,P_SB^rP_S)=\{\varepsilon x\in \mathbb{R}^r: \varepsilon\in[0,1], x\in m_S\}$ \ (see \eqref{def jnr}). If $S=\C^n$, then $W(B^1,\dots,B^r)=\Phi^{\B}(\calD)=m_{\C^n}$.
	\end{enumerate}
\end{proposition}

\begin{proof}
	\begin{enumerate}
		\item Recall that any $\mu\in \mathcal{D}_S$ is a convex combination of rank-one orthogonal projections. That means that it can be described as $\mu=\sum_{k=1}^{\dim(S)} t_k z_kz_k^*$ for $t_k\geq 0$, $\sum_{k=1}^{\dim S} t_k=1$ and $\{z_k\}_{k=1}^{\dim S}$ an orthonormal basis of $S$. Then
		$$
		\tr(\mu B^i)=\sum_{k=1}^{\dim S} t_k\tr(B^i z_kz_k^*)=\sum_{k=1}^{\dim S} t_k \langle B^i z_k, z_k\rangle, \quad\text{for } i=1,\dots, r.
		$$
		Therefore
			$$
			\begin{aligned}
				\Phi^{\B}(\mu)&=\left(\tr(B^1\mu), \dots, \tr(B^r\mu)\right)\\
				&=\sum_{k=1}^{\dim(S)}t_k\left(\langle B^1 z_k,z_k\rangle,\dots, \langle B^rz_k,z_k\rangle\right)\in m_S,
			\end{aligned}
			$$
			which proves one inclusion.
			
			Now, for every $v\in S$ with $\|v\|=1$, we have that
			$$
			(\langle B^1v,v\rangle,\dots,\langle B^r v,v\rangle)=\Phi^{\B}(vv^*), \qquad vv^*\in \mathcal{D}_S.
			$$
			Since $\mathcal{D}_S$ is convex and $\Phi^{\B}$ is linear, $\Phi^{\B}(\calD_S)$ is convex, and therefore
			$$
			m_S=\operatorname{conv}\{(\langle B^1v,v\rangle,\dots,\langle B^r v,v\rangle):v\in S,\ \|v\|=1\} \subseteq \Phi^{\B}(\calD_S),
			$$
			which proves the other inclusion.
		\item This is item \ref{mS = diag DsubS} written explicitly, by the definition of $\Phi^{\B}$.
		\item The inclusion
		$\left\{\left(\langle B^1 z_k, z_k\rangle ,
		\dots,
		\langle B^r   z_k, z_k\rangle \right) \right\}_{k=1}^{\dim(S)}\subset m_S$, for $\{z_k\}_{k=1}^{\dim(S)}$ an orthonormal basis, follows from item \ref{mS = diag DsubS}, since $z_k z_k^*\in \mathcal{D}_S$, and the fact that $m_S$ is convex. Hence the union of these sets is also included in $m_S$.
		
		The other inclusion follows from the description of every $\mu\in \mathcal{D}_S$ in terms of an orthonormal basis of $S$, as in the beginning of the proof of item \ref{mS = diag DsubS}.
		\item This statement can be proved using the same arguments as in the proof of $(4)\Rightarrow(5)$ in \cite[Proposition 1]{BV-2026}. The hypothesis $S\neq\C^n$ is needed for the inclusion $\supseteq$ with $\varepsilon<1$, which uses density matrices supported on $S^\perp$. If $S=\C^n$, then $P_S=I$ and the claim is item \ref{mS = diag DsubS} with $\calD_{\C^n}=\calD$.
		\end{enumerate}
\end{proof}

For the next result, note that 
$\overline{0}\in W(P_SB^1P_S,\dots,P_SB^rP_S)$ whenever \( S\subsetneq \C^n\), since $\rho=\frac{1}{\tr(P_{S^\perp})}P_{S^\perp}\in\mathcal D$ satisfies $\tr(P_SB^iP_S\,\rho)=0, \ \text{for all } i=1,\dots,r.$

\begin{proposition}\label{prop 0 en momento}
	With the previous notation, the following statements are equivalent for $\{0\}\neq S\neq\mathbb C^n$.
	\begin{enumerate}
		\item $\overline{0}\in m_S^\B$, 
		\item $m_S^\B=W(P_SB^1P_S,\dots,P_SB^rP_S)$,
		\item $\exists \mu\in \calD_S$ such that $\tr(\mu B^i)=0$, for $i=1,\dots, r$,
		\item $\exists \mu\in \calD_S\cap \calW^\circ_1$,
		\item $\calD_{S}\cap\calW^\circ\neq \emptyset$.
	\end{enumerate}
\end{proposition}
\begin{proof}
	All the equivalences follow by direct computation from Proposition \ref{prop momento}.
\end{proof}

The following result is proved in \cite{BV-2026} for $C^*$-subalgebras of $M_n(\C)$; the same proof applies here since it only uses orthonormal Hermitian bases of $\calW$.
	\begin{lemma}\label{phi change of basis}
		If $\mathbb{T}=\{T^k\}_{k=1}^{r}$ and $\mathbb{V}=\{V^k\}_{k=1}^{r}$ are two orthonormal bases of Hermitian matrices of $\calW$, then for any $\rho\in \calD_{S}$
		\begin{align*}
			\Phi^{\mathbb{V}}(\rho)= C_{\mathbb{T},\mathbb{V}}\Phi^{\mathbb{T}}(\rho),
		\end{align*}
		where $C_{\mathbb{T},\mathbb{V}}$ is the change of basis matrix from $\mathbb{T}$ to $\mathbb{V}$, that is
		\[(C_{\mathbb{T},\mathbb{V}})_{ij}=\langle V^i, T^j\rangle=\tr(T^jV^i), \ \text{with } i,j=1,\dots,r.
		\] 
	\end{lemma}
	\begin{proposition}[Change of basis in $\calW$]\label{prop change of basis} If $\mathbb{T}=\{T^k\}_{k=1}^{\dim(\calW)}$ and $\mathbb{V}=\{V^k\}_{k=1}^{\dim(\calW)}$ are two orthonormal bases of Hermitian matrices of $\calW$, then
		\begin{equation}\label{change of basis moment}
			m_S^{\mathbb{V}}= C_{\mathbb{T},\mathbb{V}}m_S^{\mathbb{T}}
		\end{equation}
		and
		\begin{equation}\label{change of basis jnr}
			W(P_S\mathbb{V}P_S)=C_{\mathbb{T},\mathbb{V}}W(P_S\mathbb{T}P_S),
		\end{equation}
		where $C_{\mathbb{T},\mathbb{V}}$ is the change of basis matrix.	
\end{proposition}

\begin{proof}
	Combining Lemma \ref{phi change of basis} and item 1 in Proposition \ref{prop momento} it follows that $$m_S^{\mathbb{V}}=\Phi^{\mathbb{V}}(\calD_S)= C_{\mathbb{T},\mathbb{V}}\Phi^{\mathbb{T}}(\calD_S)=C_{\mathbb{T},\mathbb{V}}m_S^{\mathbb{T}},$$
	which proves \eqref{change of basis moment}.
	If $S=\C^n$, then \eqref{change of basis jnr} is \eqref{change of basis moment}, by item 4 of Proposition \ref{prop momento}. Otherwise, \eqref{change of basis jnr} is a direct consequence of \eqref{change of basis moment} and item 4 of Proposition \ref{prop momento}, since
	\begin{align*}
		W(P_S\mathbb{V}P_S)&= \{\varepsilon x': x'\in m_S^{\mathbb{V}}, \varepsilon\in[0,1]\}\\
		&= \{\varepsilon C_{\mathbb{T},\mathbb{V}}x: x\in m_S^{\mathbb{T}}, \varepsilon\in[0,1]\}\\
		&= C_{\mathbb{T},\mathbb{V}}W(P_S\mathbb{T}P_S).
	\end{align*}	
	
\end{proof}

In the following proposition we relate the field (see \eqref{field})
\[
\mathscr{F}(\{(A_0+W_0)B^j\}_{j=1}^r, S_{\mathrm{MAX}})
\]
to the fields of $\{B^j\}_{j=1}^r$ restricted to the eigenspaces $S_{\lambda_{\max}}$ and $S_{\lambda_{\min}}$.

\begin{proposition} \label{relation sets F}
	Let $\calW\subset M_n^h(\C)$ be a subspace with $\dim_{\R}(\calW)=r$, let $\{B^j\}_{j=1}^r$ be the fixed orthonormal basis of Hermitian matrices of $\calW$, let $A_0\in M_n^h(\C)\setminus \calW$ and $W_0\in \calW$, and let $S_{\mathrm{MAX}}=S_{\mathrm{MAX}}(A_0+W_0)$ be as in \eqref{SMAX generico}. Let $\lambda_{\max}$ and $\lambda_{\min}$ be the maximum and minimum eigenvalues of $A_0+W_0$. Then:
	\begin{enumerate}
		\item If $\lambda_{\max}+\lambda_{\min}>0$ (equivalently, $S_{\mathrm{MAX}}=S_{\lambda_{\max}}$), then
		\[\mathscr{F}(\{(A_0+W_0)B^j\}_{j=1}^r, S_{\mathrm{MAX}})=\|A_0+W_0\|\,\mathscr{F}(\{B^j\}_{j=1}^r,S_{\lambda_{\max}})\subset \R^r.\]
		\item If $\lambda_{\max}+\lambda_{\min}<0$ (equivalently, $S_{\mathrm{MAX}}=S_{\lambda_{\min}}$), then
		\[\mathscr{F}(\{(A_0+W_0)B^j\}_{j=1}^r, S_{\mathrm{MAX}})=-\|A_0+W_0\|\,\mathscr{F}(\{B^j\}_{j=1}^r,S_{\lambda_{\min}})\subset \R^r.\]
		\item If $\lambda_{\max}+\lambda_{\min}=0$ (equivalently, $S_{\mathrm{MAX}}=S_{\lambda_{\max}}\oplus S_{\lambda_{\min}}$), then
		\begin{align*}
		&\mathscr{F}(\{(A_0+W_0)B^j\}_{j=1}^r, S_{\mathrm{MAX}})\\
		&\quad=\|A_0+W_0\|\left\{t\,y(u_+)-(1-t)\,z(u_-)+2\mathrm{i}\sqrt{t(1-t)}\,w(u_+,u_-)\right\},
		\end{align*}
		where the set on the right runs over $t\in[0,1]$ and unit vectors $u_+\in S_{\lambda_{\max}}$, $u_-\in S_{\lambda_{\min}}$, and
		\begin{align*}
			y(u_+)&=\left(\langle B^1u_+,u_+\rangle, \dots, \langle B^ru_+,u_+\rangle\right)\in \mathscr{F}(\{B^j\}_{j=1}^r,S_{\lambda_{\max}}),\\
			z(u_-)&=\left(\langle B^1u_-,u_-\rangle, \dots, \langle B^ru_-,u_-\rangle\right)\in \mathscr{F}(\{B^j\}_{j=1}^r,S_{\lambda_{\min}}),\\
			w(u_+,u_-)&=\left(\operatorname{Im}\langle B^1u_-,u_+\rangle, \dots, \operatorname{Im}\langle B^ru_-,u_+\rangle\right)\in\R^r.
		\end{align*}
		In this case the field need not be contained in $\R^r$, since the cross terms $\langle B^j u_-,u_+\rangle$ are in general complex.
	\end{enumerate}
\end{proposition}

\begin{proof}
	Since $A_0\notin\calW$, $A_0+W_0\neq0$, so by Proposition \ref{lema prop smax} the three conditions on $\lambda_{\max}+\lambda_{\min}$ are equivalent to the stated descriptions of $S_{\mathrm{MAX}}$. Let $v\in S_{\mathrm{MAX}}$ be a unit vector and consider
	\[x =\left( \langle (A_0+W_0)B^1v,v\rangle, \dots, \langle (A_0+W_0)B^rv,v\rangle\right),\]
	so that $x\in\mathscr{F}(\{(A_0+W_0)B^j\}_{j=1}^r, S_{\mathrm{MAX}})$, and every element of this field is of this form.
	Since $B^j$ and $A_0+W_0$ are Hermitian, $\langle (A_0+W_0)B^jv,v\rangle=\langle B^jv,(A_0+W_0)v\rangle$ for every $1\leq j\leq r$.
	\begin{enumerate}
		\item If $S_{\mathrm{MAX}}=S_{\lambda_{\max}}$, then $(A_0+W_0)v=\|A_0+W_0\|v$ for every $v\in S_{\mathrm{MAX}}$, so
		\[x= \|A_0+W_0\| \left( \langle B^1v,v\rangle, \dots, \langle B^rv,v\rangle\right)\in \|A_0+W_0\|\mathscr{F}(\{B^j\}_{j=1}^r,S_{\lambda_{\max}}).\]
		Reversing the same steps gives the opposite inclusion.
		\item If $S_{\mathrm{MAX}}=S_{\lambda_{\min}}$, the same argument, with $(A_0+W_0)v=-\|A_0+W_0\|v$, gives
		\[x\in -\|A_0+W_0\|\mathscr{F}(\{B^j\}_{j=1}^r,S_{\lambda_{\min}}),\]
		and again the opposite inclusion follows by reversing the steps.
		\item If $\lambda_{\max}+\lambda_{\min}=0$, then $\|A_0+W_0\|=\lambda_{\max}$. Write $v=v_++v_-$ with $v_+\in S_{\lambda_{\max}}$ and $v_-\in S_{\lambda_{\min}}$. Since the eigenspaces are orthogonal, $t=\|v_+\|^2\in[0,1]$ and $\|v_-\|^2=1-t$, so that $v_+=\sqrt{t}\,u_+$ and $v_-=\sqrt{1-t}\,u_-$ for some unit vectors $u_+\in S_{\lambda_{\max}}$ and $u_-\in S_{\lambda_{\min}}$ (if $t\in\{0,1\}$, the vector multiplied by $0$ is an arbitrary unit vector of the corresponding eigenspace). Since $(A_0+W_0)v=\lambda_{\max}(v_+-v_-)$, $B^j$ is Hermitian and $\langle B^jv_+,v_-\rangle=\overline{\langle B^jv_-,v_+\rangle}$, we have that
		\begin{align*}
			\langle (A_0+W_0)B^jv,v\rangle &=\lambda_{\max}\langle B^j(v_++v_-),v_+-v_-\rangle\\
			&=\lambda_{\max}\big[\langle B^jv_+,v_+\rangle-\langle B^jv_-,v_-\rangle\\
			&\qquad\qquad+\langle B^jv_-,v_+\rangle-\langle B^jv_+,v_-\rangle\big]\\
			&=\lambda_{\max}\big[t\langle B^ju_+,u_+\rangle-(1-t)\langle B^ju_-,u_-\rangle\\
			&\qquad\qquad+2\mathrm{i}\sqrt{t(1-t)}\,\operatorname{Im}\langle B^ju_-,u_+\rangle\big]
		\end{align*}
		for every $1\leq j\leq r$. Hence
		\[x=\lambda_{\max}\left[t\,y(u_+)-(1-t)\,z(u_-)+2\mathrm{i}\sqrt{t(1-t)}\,w(u_+,u_-)\right],\]
		which proves the inclusion $\subseteq$. Conversely, given $t\in[0,1]$ and unit vectors $u_\pm\in S_{\lambda_{\max/\min}}$, the vector $v=\sqrt t\,u_++\sqrt{1-t}\,u_-$ is a unit vector of $S_{\mathrm{MAX}}$, and the same computation shows that the corresponding element of the right-hand side belongs to the field.
	\end{enumerate}
\end{proof}

\begin{corollary}\label{cor field y momento}
	Under the hypotheses of Proposition \ref{relation sets F}, let $\lambda_{\max}$ and $\lambda_{\min}$ be the maximum and minimum eigenvalues of $A_0+W_0$. Recall that for every subspace $S\subset \C^n$, $\operatorname{conv}\mathscr{F}(\{B^j\}_{j=1}^r,S)=m_S$ (Remark \ref{obs sobre momento}). Hence,
	\begin{equation}\label{mom=conv F lambdas}
		\begin{split}
			m_{S_{\lambda_{\max}}}&=\operatorname{conv}\mathscr{F}(\{B^j\}_{j=1}^r,S_{\lambda_{\max}}),\\  m_{S_{\lambda_{\min}}}&=\operatorname{conv}\mathscr{F}(\{B^j\}_{j=1}^r,S_{\lambda_{\min}}),
		\end{split}
	\end{equation}
	and, taking real parts componentwise in $\C^r$,
	\begin{equation}\label{momento=conv F}
		\begin{split}
			&\operatorname{Re}\left(\operatorname{conv}\mathscr{F}(\{(A_0+W_0)B^j\}_{j=1}^r, S_{\mathrm{MAX}})\right)\\
			&\qquad=\|A_0+W_0\|\cdot\left\lbrace \begin{array}{ll}
				m_{S_{\mathrm{MAX}}} & \text{if } \lambda_{\max}+\lambda_{\min}>0,\\
				-m_{S_{\mathrm{MAX}}} & \text{if } \lambda_{\max}+\lambda_{\min}<0,\\
				\operatorname{conv}\left(m_{S_{\lambda_{\max}}}\cup -m_{S_{\lambda_{\min}}}\right) & \text{if } \lambda_{\max}+\lambda_{\min}=0.
			\end{array}\right.
		\end{split}
	\end{equation}
\end{corollary}
\begin{proof}
	Write
	\begin{gather*}
		\mathscr{F}_{\mathrm{MAX}}=\mathscr{F}(\{(A_0+W_0)B^j\}_{j=1}^r, S_{\mathrm{MAX}}),\\
		\mathscr{F}_+=\mathscr{F}(\{B^j\}_{j=1}^r,S_{\lambda_{\max}}),\qquad \mathscr{F}_-=\mathscr{F}(\{B^j\}_{j=1}^r,S_{\lambda_{\min}}).
	\end{gather*}
	Since $\operatorname{Re}$ is $\R$-linear, $\operatorname{Re}(\operatorname{conv}X)=\operatorname{conv}(\operatorname{Re}X)$ for every $X\subseteq\C^r$. If $\lambda_{\max}+\lambda_{\min}>0$, then $S_{\mathrm{MAX}}=S_{\lambda_{\max}}$, and by item 1 of Proposition \ref{relation sets F} we have $\mathscr{F}_{\mathrm{MAX}}\subset\R^r$ and
	\[\operatorname{conv}\mathscr{F}_{\mathrm{MAX}}=\|A_0+W_0\|\operatorname{conv}\mathscr{F}_+=\|A_0+W_0\|\,m_{S_{\mathrm{MAX}}}.\]
	The case $\lambda_{\max}+\lambda_{\min}<0$ follows in the same way from item 2.
	If $\lambda_{\max}+\lambda_{\min}=0$, the computation in the proof of item 3 holds for every $t\in[0,1]$ and every pair of unit vectors $u_+\in S_{\lambda_{\max}}$, $u_-\in S_{\lambda_{\min}}$, taking $v=\sqrt{t}\,u_++\sqrt{1-t}\,u_-$. Therefore
	\[\operatorname{Re}\mathscr{F}_{\mathrm{MAX}}=\|A_0+W_0\|\left\{t\,y-(1-t)\,z:\ t\in[0,1],\ y\in\mathscr{F}_+,\ z\in\mathscr{F}_-\right\}.\]
	The set between braces contains $\mathscr{F}_+$ (take $t=1$) and $-\mathscr{F}_-$ (take $t=0$), and it is contained in $\operatorname{conv}(\mathscr{F}_+\cup-\mathscr{F}_-)$. Hence its convex hull is $\operatorname{conv}(\mathscr{F}_+\cup-\mathscr{F}_-)=\operatorname{conv}(m_{S_{\lambda_{\max}}}\cup -m_{S_{\lambda_{\min}}})$, by \eqref{mom=conv F lambdas}.
\end{proof}

\begin{remark}\label{rem field y minimalidad}
By Corollary \ref{formula for subdif composition norm} below, the right-hand side of \eqref{momento=conv F} is $\|A_0+W_0\|$ times $\partial(\|A(\cdot)\|)(x)$, where $x\in\R^r$ is such that $A(x)=A_0+W_0$. Since $\|A_0+W_0\|>0$, Theorem \ref{teo minimality merged} shows that $A_0+W_0$ is $\calW$-minimal if and only if
\[
\overline{0}\in\operatorname{Re}\left(\operatorname{conv}\mathscr{F}(\{(A_0+W_0)B^j\}_{j=1}^r, S_{\mathrm{MAX}})\right).
\]
\end{remark}

\section{Subdifferentials}\label{sec:subdiff}

First, we introduce the definition of the subdifferential of a convex function over a Banach space. 
\begin{definition}\label{def subdiferencial}
	Let $\mathcal{X}$ be a Banach space with dual $\mathcal{X}^*$, and let $f:\mathcal{X}\to \R$ be a convex function. The subdifferential of $f$ at $x\in\mathcal{X}$ is
	\begin{equation}\label{subdif gral banach}
		\partial f(x)=\{v\in\mathcal{X}^*: f(y)-f(x)\geq \operatorname{Re} v(y-x), \forall\ y\in \mathcal{X}\}.
	\end{equation}
\end{definition}
In particular for $f:M_n(\C)\to \R$, $f(A)=\|A\|$, 
\[	\partial \|A\|=\{G\in M_n(\C): \|B\|-\|A\|\geq \operatorname{Re} \tr(G(B-A)), \forall\ B\in M_n(\C)\},\]
and it is well known that this definition is equivalent to
\begin{equation}\label{subdif norma def}
\partial \|A\|=\{G\in M_n(\C): \|A\|=\operatorname{Re} \tr(G^*A), \ \|G\|_1\leq 1\}.
\end{equation}
Moreover, by \cite[Theorem 4]{watson}, the subdifferential of the spectral norm of any $A\in M_n(\C)$ has the equivalent formulation
\begin{equation} \label{watson}
	\partial\|A\|=\operatorname{conv}\{uv^*:\ Av=\|A\|u,\ \|v\|=\|u\|=1\}.
\end{equation}
For a set $\mathcal A\subseteq M_n(\C)$, we write $\operatorname{Re}(\mathcal A)=\{\operatorname{Re}(X): X\in\mathcal A\}$, where $\operatorname{Re}(X)=\frac{X+X^*}{2}$. We obtain the following result for Hermitian matrices.
\begin{lemma}\label{lema subdif hermitiana}
		Let $A\in M_n^h(\C)$ and let $G=\operatorname{Re}(G)+\mathrm{i}\operatorname{Im}(G)$ be
		the Cartesian decomposition of a matrix $G\in M_n(\C)$, that is,
		$\operatorname{Re}(G),\operatorname{Im}(G)\in M_n^h(\C)$. Then, with respect to
		the real inner product $\langle X,Y\rangle=\operatorname{Re}\tr(X^*Y)$,
		\begin{equation}\label{eq subdif}
			G\in\partial\|A\|\ \Leftrightarrow\ \tr(\operatorname{Re}(G)A)=\|A\|
			\ \text{ and }\ \|G\|_1\leq 1.
		\end{equation}
		Moreover, if $\partial_h\|A\|$ denotes the subdifferential of the restriction
		of $\|\cdot\|$ to $M_n^h(\C)$, then
		\[
		\partial_h\|A\|=\partial\|A\|\cap M_n^h(\C)=\operatorname{Re}(\partial\|A\|)
		\]
		and
		\begin{equation}\label{eq subdif herm}
			\partial_h\|A\|=\operatorname{conv}\{s\,vv^*:\ s\in\{1,-1\},\ \|v\|=1,\
			Av=s\|A\|v\}.
		\end{equation}
	\end{lemma}
	\begin{proof}
		By \eqref{subdif norma def}, $G\in\partial\|A\|$ if and only if
		$\operatorname{Re}\tr(G^*A)=\|A\|$ and $\|G\|_1\leq 1$. Since
		$G^*=\operatorname{Re}(G)-\mathrm{i}\operatorname{Im}(G)$ and the trace of a product of
		two Hermitian matrices is real,
		$\operatorname{Re}\tr(G^*A)=\tr(\operatorname{Re}(G)A)$, which proves
		\eqref{eq subdif}.
		\\
		Clearly $\partial\|A\|\cap M_n^h(\C)\subseteq\partial_h\|A\|$. Conversely, let
		$G\in\partial_h\|A\|$ and $B\in M_n(\C)$. Since $G$ is Hermitian,
		$\langle G,B\rangle=\tr(G\operatorname{Re}(B))$. Applying the subgradient
		inequality to the Hermitian matrix $\operatorname{Re}(B)$ and using
		$\|\operatorname{Re}(B)\|\le\|B\|$, we get
		\[
		\|A\|+\langle G,B-A\rangle=\|A\|+\langle G,\operatorname{Re}(B)-A\rangle
		\le\|\operatorname{Re}(B)\|\le\|B\|,
		\]
		so $G\in\partial\|A\|$. Hence $\partial_h\|A\|=\partial\|A\|\cap M_n^h(\C)$.
		If $G\in\partial\|A\|$, then $G^*\in\partial\|A\|$ by \eqref{eq subdif}, since
		$\|G^*\|_1=\|G\|_1$ and $\operatorname{Re}(G^*)=\operatorname{Re}(G)$. Then, by
		convexity, $\operatorname{Re}(G)=\frac{G+G^*}{2}\in\partial\|A\|\cap M_n^h(\C)$.
		As $\operatorname{Re}(G)=G$ for Hermitian $G$, it follows that
		$\partial_h\|A\|=\operatorname{Re}(\partial\|A\|)$.
		\\
		If $A=0$, then $\partial_h\|0\|=\{G\in M_n^h(\C):\|G\|_1\le1\}$ by
		\eqref{eq subdif}, and every unit vector is admissible in
		\eqref{eq subdif herm} with both signs. Writing $G=\sum_i\mu_iw_iw_i^*$, with
		$\{w_i\}$ orthonormal, we have
		$G=\sum_i|\mu_i|\,(\operatorname{sgn}(\mu_i)\,w_iw_i^*)
		+\tfrac d2\,vv^*+\tfrac d2\,(-vv^*)$, where $d=1-\|G\|_1\ge0$ and $v$ is any
		unit vector. This proves that $G\in\operatorname{conv}\{\pm vv^*:\ \|v\|=1\}$, hence
		$\partial_h\|0\|\subseteq\operatorname{conv}\{\pm vv^*:\ \|v\|=1\}$. The
		reverse inclusion holds because $\partial_h\|0\|$ is convex and contains
		every $\pm vv^*$ with $\|v\|=1$. Therefore \eqref{eq subdif herm} holds in
		this case.
		\\
		Assume now $A\neq0$. By Watson's formula \eqref{watson}, $\partial\|A\|$ is the convex hull of the matrices $uv^*$ with $\|u\|=\|v\|=1$ and $Av=\|A\|u$.
		For such $v$ we have $v\in M_A$, so by Proposition~\ref{lema prop smax} we can
		write $v=v_++v_-$ with $v_\pm\in S_{\pm\|A\|}$, and then
		$u=\frac{Av}{\|A\|}=v_+-v_-$. Since $\operatorname{Re}$ is linear,
		$\partial_h\|A\|=\operatorname{Re}(\partial\|A\|)$ is the convex hull of the
		matrices
		\[
		\operatorname{Re}\bigl((v_+-v_-)(v_++v_-)^*\bigr)=v_+v_+^*-v_-v_-^*
		=\|v_+\|^2\,w_+w_+^*+\|v_-\|^2\,(-w_-w_-^*),
		\]
		where $w_\pm=v_\pm/\|v_\pm\|$ (omit a term if $v_\pm=0$) and
		$\|v_+\|^2+\|v_-\|^2=\|v\|^2=1$. Each of these is a convex combination of
		matrices $s\,ww^*$ with $\|w\|=1$ and $Aw=s\|A\|w$. Conversely, each such
		$s\,ww^*$ is of this form (take $v=w$, so that $v_-=0$ if $s=1$ and $v_+=0$ if
		$s=-1$). This proves \eqref{eq subdif herm}.
\end{proof}

\begin{lemma} \label{lema of norm subdif} 
	Let $A\in M_n^h(\C)\setminus\{0\}$, and let the sets $\calT$, $\calT_+$ and $\calT_-$ be defined as
	\begin{equation}\label{def T}
		\begin{array}{ll}
			\calT&= \{uv^*:\ Av=\|A\|u, \|v\|=\|u\|=1\},\\
			\calT_+&=\{uu^*:\ Au=\|A\|u, \|u\|=1\},\\
			\calT_-&=\{-uu^*:\ Au=-\|A\|u, \|u\|=1\}.
		\end{array}
	\end{equation}
	Then,	 $\calT_+\cap \calT_-=\emptyset$ and
	
	\begin{equation}\label{subdif}
		\partial\|A\|={\rm conv}(\calT)={\rm conv}\left\{\tfrac{1}{\|A\|}Av\,v^*:\ v\in S_{\mathrm{MAX}}(A),\ \|v\|=1\right\}.
	\end{equation}
	Additionally,
	\begin{equation}\label{subdif no balanceado}
		\partial\|A\|=\left\lbrace \begin{array}{lll}
			{\rm conv}(\calT_+)&\text{if}& \lambda_{\max}=\|A\|>-\lambda_{\min},\\
			{\rm conv}(\calT_-)&\text{if}& -\lambda_{\min}=\|A\|>\lambda_{\max},\\
		\end{array}\right. 					
	\end{equation}
	and if $\lambda_{\max}=-\lambda_{\min}=\|A\|$
	\begin{equation} \label{subdif espectro balanceado}
		\partial\|A\|=\operatorname{conv}\{(v_+-v_-)(v_++v_-)^*: v_+\in S_{\|A\|}, v_-\in S_{-\|A\|},\|v_+\|^2+\|v_-\|^2=1 \}.				
	\end{equation}
	Finally,
	\begin{equation}\label{subdif 2}
		\partial_h\|A\|=\operatorname{Re}(\partial\|A\|)={\rm conv}(\calT_+\cup \calT_-).
	\end{equation}
\end{lemma}

\begin{proof}
	First observe that if $v$ is a nonzero vector in $S_{\|A\|}$ (resp. $S_{-\|A\|}$), then $v\notin S_{-\|A\|}$ ($v\notin S_{\|A\|}$), so $vv^*\in \calT_+$ ($vv^*\in \calT_-$) but $vv^*\notin \calT_-$ ($vv^*\notin \calT_+$). Therefore, $\calT_+\cap \calT_-=\emptyset$. 
	
	By \eqref{watson}, $\partial\|A\|={\rm conv}(\calT)$. Moreover, if $uv^*\in\calT$, then $Av=\|A\|u$ with $\|v\|=\|u\|=1$, so $\|Av\|=\|A\|$, $v\in M_A\subseteq S_{\mathrm{MAX}}(A)$ by item 4 of Proposition \ref{lema prop smax}, and $u=\frac{1}{\|A\|}Av$. Conversely, if $v\in S_{\mathrm{MAX}}(A)$ and $\|v\|=1$, then $\|Av\|=\|A\|$ by the same item, so $u=\frac{1}{\|A\|}Av$ is a unit vector with $Av=\|A\|u$, and $uv^*\in\calT$. This proves \eqref{subdif}. Now fix $uv^*\in\calT$. We consider the following cases: 
	\begin{itemize}
		\item If $\lambda_{\max}=\|A\|>-\lambda_{\min}$: here $S_{\mathrm{MAX}}(A)=S_{\|A\|}$ and it follows that $u=\frac{1}{\|A\|}Av=v$. Moreover, $uv^*=vv^*\in \calT_+$, then $\calT=\calT_+$ and $\partial\|A\|={\rm conv}(\calT_+)$.
		\item If $-\lambda_{\min}=\|A\|>\lambda_{\max}$: in this case $S_{\mathrm{MAX}}(A)=S_{-\|A\|}$ and it follows that $u=\frac{1}{\|A\|}Av=-v$. Moreover, $uv^*=-vv^*\in \calT_-$, then $\calT=\calT_-$ and $\partial\|A\|={\rm conv}(\calT_-)$.
		\item If $\lambda_{\max}=\|A\|=-\lambda_{\min}$: $v\in S_{\mathrm{MAX}}(A)=S_{\|A\|}\oplus S_{-\|A\|}$, then there exist $v_+\in S_{\|A\|}$ and $v_-\in S_{-\|A\|}$ such that $v=v_++v_-$ and $\|v\|^2=\|v_+\|^2+\|v_-\|^2=1$. By item 5 in Proposition \ref{lema prop smax}, $u=\frac{1}{\|A\|}Av=v_+-v_-$. Moreover,  
		\[uv^*=(v_+-v_-)(v_++v_-)^*,\]
		which proves \eqref{subdif espectro balanceado}.
		
	\end{itemize}
	Finally, by \eqref{def T}, $\calT_+\cup\calT_-=\{s\,vv^*:\ s\in\{1,-1\},\ \|v\|=1,\ Av=s\|A\|v\}$, so \eqref{subdif 2} is a restatement of Lemma \ref{lema subdif hermitiana}.
\end{proof}

We recall the following chain rule, which we use to obtain a closed form for the subdifferential of the norm of any Hermitian $A(x)$.
\begin{lemma}[Chain rule \cite{clarke1, grover26}]\label{regla de la cadena}
	Let $\calX$ and $\calY$ be Banach spaces over the field $\R$ (or $\C$). Let $S:\calX\to \calY$ be a bounded
	linear map, and let $L:\calX\to \calY$ be the affine map defined as $L(x) = S(x) + y_0$, where $y_0\in\calY$. If $g:\calY\to \R$ is a continuous convex function, then
	\[\partial(g\circ L)(a)=S^*\partial g(L(a)), \text{ for all } a\in \calX,\]
	where $S^*$ denotes the adjoint of the operator $S$, either real or complex, depending on whether $\calX$ and $\calY$ 	are real or complex Banach spaces.
\end{lemma}

\begin{theorem} \label{teo subdif norma}
	For $A(x)=A_0+\sum_{k=1}^{r}x_kB^k$, with $A_0\in M_n^{h}(\C)$ and $\{B^k\}_{k=1}^{r}$ the fixed orthonormal basis of $\calW$, $x\in \R^{r}$, let $\lambda_{\min}(A(x))$ and $\lambda_{\max}(A(x))$ be the minimum and maximum eigenvalues of $A(x)$, respectively. Then, the subdifferential of the spectral norm of $A(x)$ is
	\begin{equation}\label{watsonhermit}
		\begin{split}
			\partial(\|A(\cdot)\|)(x)&=\Phi\left(\operatorname{Re}(\partial\|A(x)\|)\right)=\Phi\left({\rm conv}(\calT_+\cup\calT_-)\right)\\
			&=\left\lbrace \begin{array}{lll}
				\Phi\left( {\rm conv} \left( \calT_+\right)\right) &\text{ if }& \lambda_{\max}+\lambda_{\min}>0\\			
				\Phi\left( {\rm conv} \left( \calT_-\right)\right) &\text{ if }& \lambda_{\max}+\lambda_{\min}<0\\			
				\Phi\left( {\rm conv} \left( \calT_+\cup \calT_-\right)\right) &\text{ if }& \lambda_{\max}+\lambda_{\min}=0
			\end{array}\right.  ,			
		\end{split}
	\end{equation}
	where $\lambda_{\max}=\lambda_{\max}(A(x))$, $\lambda_{\min}=\lambda_{\min}(A(x))$, $\Phi=\Phi^\B$ is the map defined in \eqref{defi phi general}, and $\calT_+$, $\calT_-$ are the sets defined in \eqref{def T} for $A=A(x)$.

\end{theorem}

\begin{proof}
	Apply Lemma \ref{regla de la cadena} with $\calX=\R^{r}$, $\calY=M_n^h(\C)$ endowed with the real inner product $\langle X,Y\rangle=\tr(XY)$, $S=W$, where $W(x)=\sum_{k=1}^{r}x_kB^k$, $L(x)=A(x)=A_0+W(x)$, and $g:M_n^h(\C)\to\R$, $g(R)=\|R\|$. Then $g\circ L=\|A(\cdot)\|$ and $\partial g(A(x))=\partial_h\|A(x)\|$ is the subdifferential of the restriction of the norm to $M_n^h(\C)$. For $x\in\R^r$ and $R\in M_n^h(\C)$,
	\[
	\langle W(x),R\rangle=\sum_{k=1}^{r}x_k\tr(B^kR)=\langle x,\Phi(R)\rangle_{\R^r},
	\]
	so $W^*=\Phi$ and
	\begin{equation}\label{RDC}
		\partial(\|A(\cdot)\|)(x)=\Phi\left(\partial_h\|A(x)\|\right).
	\end{equation}
	By Lemma \ref{lema subdif hermitiana} we have $\partial_h\|A(x)\|=\operatorname{Re}(\partial\|A(x)\|)$, and by \eqref{subdif 2} in Lemma \ref{lema of norm subdif} we have $\operatorname{Re}(\partial\|A(x)\|)={\rm conv}(\calT_+\cup\calT_-)$. Finally, if $\lambda_{\max}+\lambda_{\min}>0$ (resp. $<0$) for $A(x)$, then $\calT_-=\emptyset$ (resp. $\calT_+=\emptyset$), which gives the three cases in \eqref{watsonhermit}.
\end{proof}
Formulas for the subdifferentials of $\lambda_{\max}$ and $\lambda_{\min}$ at $A$ were given in \cite{overton2, BV-SDP, BV-2026}. These are
\begin{equation} \label{subdif eig matrices}
	\begin{split}
		\partial\lambda_{\max}(A)&=\operatorname{conv}\{qq^*: Aq=\lambda_{\max}(A)q \text{ and }\|q\|=1\}\\
		\partial\lambda_{\min}(A)&=-\partial\lambda_{\max}(-A)=-\operatorname{conv}\{qq^*: Aq=\lambda_{\min}(A)q \text{ and }\|q\|=1\}.
	\end{split}
\end{equation}
If $S_{\lambda_{\max}}$ and $S_{\lambda_{\min}}$ are the eigenspaces of $\lambda_{\max}(A(x))$ and $\lambda_{\min}(A(x))$, respectively, in \cite[Theorem 3 and Proposition 3]{BV-2026} we showed that
\begin{equation}\label{subdiff max y min = momento}
	\partial(\lambda_{\max}\circ A)(x)= m_{S_{\lambda_{\max}}}\ \text{and}\	\partial(\lambda_{\min}\circ A)(x)= -m_{S_{\lambda_{\min}}}.
\end{equation}
The next result follows directly from the formulas \eqref{watsonhermit} and \eqref{subdiff max y min = momento}.
\begin{corollary} \label{formula for subdif composition norm}
	Under the same hypotheses as in Theorem \ref{teo subdif norma}, the following statements hold for $A(x)$:
	\begin{enumerate}
		\item If $\|A(x)\|=\lambda_{\max}(A(x))> -\lambda_{\min}(A(x))$, then $$\partial(\|A(\cdot)\|)(x)=\partial(\lambda_{\max}\circ A)(x)=m_{S_{\lambda_{\max}}}=m_{S_{\mathrm{MAX}}(A(x))}.$$
		\item If $\|A(x)\|=-\lambda_{\min}(A(x))> \lambda_{\max}(A(x))$, then $$\partial(\|A(\cdot)\|)(x)=\partial(\lambda_{\min}\circ A)(x)=-m_{S_{\lambda_{\min}}}=-m_{S_{\mathrm{MAX}}(A(x))}.$$
		\item If $\|A(x)\|=-\lambda_{\min}(A(x))=\lambda_{\max}(A(x))$, then 
		\begin{align*}
		\partial(\|A(\cdot)\|)(x)&={\rm conv}\left(\partial(\lambda_{\max}\circ A)(x)\cup\partial(\lambda_{\min}\circ A)(x) \right)\\
		&= {\rm conv}\left(m_{S_{\lambda_{\max}}}\cup-m_{S_{\lambda_{\min}}} \right).
		\end{align*}
	\end{enumerate}
	In all three cases, with the convention $m_{\{0\}}=\emptyset$ and $S_{\pm\|A(x)\|}=\{v\in\C^n:\ A(x)v=\pm\|A(x)\|v\}$,
	\[\partial(\|A(\cdot)\|)(x)={\rm conv}\left(m_{S_{\|A(x)\|}}\cup-m_{S_{-\|A(x)\|}}\right).
	\]
\end{corollary}

\section{Characterizations of $\calW$-minimal matrices} \label{sec:characterization}
In this section we prove several characterizations of $\calW$-minimal matrices. Throughout the section we drop the superscript and write $\Phi$ for the map $\Phi^{\B}$ of \eqref{defi phi general}, defined on the whole of $M_n^h(\C)$, as well as for its restriction \eqref{defi phi} to the set $\calD_S$ of density matrices supported in a subspace $S$. 
\begin{theorem}\label{teo AX modulo}
	Let $\calW\subseteq M_n^h(\C)$ be a subspace. Then $A\in M_n^h(\C)$ is $\calW$-minimal if and only if there exists $X \in \calW^{\circ}_1$ such that
	\[AX =\|A\|\, |X|,\]
	where $|X|$ denotes the modulus $(X^*X)^{1/2}$.
\end{theorem}
\begin{proof}
	This is \cite[Theorem 2.6]{zhang-jiang-han-matrices}, and the proof given there
	remains valid not only for $C^*$-subalgebras of $M_n(\C)$, but also for arbitrary real subspaces of $M_n^h(\C)$, since the computations do not use the product of $\calW$.	
\end{proof}

When $I\in\calW$, we derive the following characterization involving the intersection of the moments of the eigenspaces of the maximum and minimum eigenvalues. Observe that if $I\in \calW$ and $A\neq0$ is $\calW$-minimal, then $\{\|A\|,-\|A\|\}\subseteq \sigma (A)$ (see also \cite{alrv}).

\begin{theorem}\label{teo mini equiv general}
	Let $A\in M_n^h(\C)$, $A\neq 0$, $\calW\subseteq M_n^h(\C)$ be a real subspace
	with $I\in\calW$, and let $E_+$ (respectively $E_-$) be the spectral projection of $A$ corresponding to
	the eigenvalue $\lambda_{\max}=\lambda_{\max}(A)$ (respectively $\lambda_{\min}=\lambda_{\min}(A)$). Suppose that $\lambda_{\max}+\lambda_{\min}=0$, and let $\B=\{B^i\}_{i=1}^{r}$ be a fixed orthonormal basis of Hermitian matrices of $\calW$, with respect to which the moments $m_S=m_S^{\B}$ of Definition \ref{momento} are computed. By \eqref{change of basis moment} the conditions below do not depend on this choice, since a change of basis acts on the moments through an invertible matrix.
	The following statements are equivalent:
	\begin{enumerate}
		\item $A$ is $\calW$-minimal.
		\item There exists $X\in\calW_1^{\circ}$ such that
		\begin{equation}\label{X witness}
			\begin{split}
				&X=E_+XE_++E_-XE_-,\quad \text{with } E_+XE_+\geq 0, \quad E_-XE_-\leq 0 \\
				&\tr(AX)=\|A\|\,\|X\|_1.
			\end{split}
		\end{equation}
		\item $m_{S_{\lambda_{\max}}}\cap m_{S_{\lambda_{\min}}}\neq\emptyset$.
	\item For every $W\in\calW$,
	\[\min_{y\in S_{\lambda_{\max}},\ \|y\|=1}\langle Wy,y\rangle\leq\max_{z\in S_{\lambda_{\min}},\ \|z\|=1}\langle Wz,z\rangle.\]
	\end{enumerate}
\end{theorem}
\begin{proof}
	We prove the implications $(1)\Rightarrow(2)\Rightarrow(3)\Rightarrow(4)$ and $(4)\Rightarrow(3)\Rightarrow(2)\Rightarrow(1)$.\\
	$(1)\Rightarrow(2)$ 
	Assume that $A$ is $\calW$-minimal, and suppose without loss of generality that
	$\|A\|=1$. 	The classical characterization of best approximation in a normed space,
	see \cite[Theorem 1.1]{singer}, provides a norm one functional that
	annihilates $\calW$ and attains $\|A\|$ at $A$. Under the trace duality
	$(M_n(\C),\|\cdot\|)^{*}=(M_n(\C),\|\cdot\|_1)$ this functional is given
	by a Hermitian matrix $X$ with
	\begin{equation}\label{eq X general}
		X\in\calW_1^{\circ},\qquad \tr(AX)=\operatorname{dist}(A,\calW)=1.
	\end{equation}	
	It remains to prove the support conditions, and for this the argument of
	\cite[Theorem 2]{alrv} applies verbatim, so we include it here for the sake of completeness. Put $E_1=E_{+}$, $E_2=E_{-}$
	and $E_3=I-E_1-E_2$, so that $E_3$ is the spectral projection of $A$
	corresponding to the open interval $(-1,1)$; then $E_3A=AE_3$,
	$\|AE_3\|<1$ and $A=E_1-E_2+AE_3$. Write $X=(X_{k,l})_{k,l=1}^{3}$ and
	$A=(A_{k,l})_{k,l=1}^{3}$ in the orthogonal decomposition
	$\C^n=\ran(E_1)\oplus\ran(E_2)\oplus\ran(E_3)$. The pinching inequality of
	Chandler Davis, see \cite[IV.52]{bhatia}, gives, writing
	$\widehat X=\begin{pmatrix}X_{1,1}&X_{1,2}\\X_{2,1}&X_{2,2}\end{pmatrix}$
	for the corner of $X$ along $\ran(E_1)\oplus\ran(E_2)$,
	\[\|\widehat X\|_1+\|X_{3,3}\|_1\leq\|X\|_1
	\qquad\text{and}\qquad
	\|X_{1,1}\|_1+\|X_{2,2}\|_1\leq\|\widehat X\|_1 .\]
	Together with $\|A_{3,3}\|<1$ and the identity $\tr(AX)=1=\|X\|_1$ these
	inequalities force $X_{3,3}=0$ as well as
	$\tr(X_{1,1})=\|X_{1,1}\|_1$ and $\tr(-X_{2,2})=\|X_{2,2}\|_1$, whence
	$X_{1,1}\geq 0$ and $X_{2,2}\leq 0$; applying the same argument to $AX$
	and then to $\widehat X$ yields $X_{k,l}=0$ for
	$k\neq l$. Therefore $X=X_{1,1}\oplus X_{2,2}\oplus 0$ with
	$X_{+}=X_{1,1}=E_{+}XE_{+}$ and $X_{-}=-X_{2,2}=-E_{-}XE_{-}$, which is
	\eqref{X witness}.\\	
	$(2)\Rightarrow(3)$ Let $X$ be as in item (2), and put $X_{+}=E_{+}XE_{+}\geq 0$ and $X_{-}=-E_{-}XE_{-}\geq 0$, so that $X=X_{+}-X_{-}$. Since $I\in\calW$ and
	$X\in\calW_1^{\circ}$ we have $\tr(X)=0$, hence
	$t=\tr(X_{+})=\tr(X_{-})$, and $t>0$ because $X\neq 0$. Then,
	\begin{equation*}
		\begin{split}
		\rho_1&=t^{-1}X_{+}=t^{-1}E_{+}XE_{+}\in\calD_{S_{\lambda_{\max}}}, \\ 
		\rho_2&=t^{-1}X_{-}=-t^{-1}E_{-}XE_{-}\in\calD_{S_{\lambda_{\min}}}  
		\end{split}
	\end{equation*}
	and $\tr(XW)=0$ for every
	$W\in\calW$.
	Therefore 
	$\tr(\rho_1W)=\tr(\rho_2W)$ for every  $W\in\calW$, and equivalently $\tr(\rho_1B^{i})=\tr(\rho_2B^i),\ 1\leq i\leq r$,
	for any Hermitian orthonormal basis $\{B^i\}_{i=1}^r$  of $\calW$.	Therefore, $\Phi(\rho_1)=\Phi(\rho_2)\in m_{S_{\lambda_{\max}}}\cap m_{S_{\lambda_{\min}}} $.
	\\	
	$(3)\Rightarrow(4)$ Let $\rho\in \calD_{S_{\lambda_{\max}}}$ and $\sigma\in \calD_{S_{\lambda_{\min}}}$ be such that $\Phi(\rho)=\Phi(\sigma)$, and let $W\in\calW$. The matrix $\rho$ is a convex combination of rank-one projections onto unit vectors of $S_{\lambda_{\max}}$, and similarly for $\sigma$ with vectors of $S_{\lambda_{\min}}$, so that
	\[\min_{y\in S_{\lambda_{\max}},\|y\|=1}\langle Wy,y\rangle\leq\tr(\rho W)=\tr(\sigma W)\leq\max_{z\in S_{\lambda_{\min}},\|z\|=1}\langle Wz,z\rangle.\]
	$(4)\Rightarrow(3)$ The sets $\calD_{S_{\lambda_{\max}}}$ and $\calD_{S_{\lambda_{\min}}}$ are
	compact and convex, and so are their images $m_{S_{\lambda_{\max}}}$ and
	$m_{S_{\lambda_{\min}}}$ under the linear map
	$\Phi$, which take values in the finite
	dimensional space $\R^r$. If these images were disjoint they could be strictly
	separated by a linear functional on $\R^r$; replacing this functional by its opposite if necessary,
	this means that there exists an element
	$x=(x_1,x_2,\dots,x_r)\in \R^r$ such that
	\[\min_{Y\in\calD_{S_{\lambda_{\max}}}}\langle \Phi(Y),x\rangle>\max_{Z\in\calD_{S_{\lambda_{\min}}}}\langle \Phi(Z),x\rangle.\]
	Taking $W=W(x)=\sum_{i=1}^{r}x_iB^i$, we get
	\[
	\min_{Y\in\calD_{S_{\lambda_{\max}}}}\tr(Y W)>\max_{Z\in\calD_{S_{\lambda_{\min}}}}\tr(ZW).
	\]
	By its spectral decomposition, every $Y\in\calD_{S_{\lambda_{\max}}}$ is a convex combination of rank-one
	projections $yy^*$ onto unit vectors $y\in S_{\lambda_{\max}}$, and these projections belong to
	$\calD_{S_{\lambda_{\max}}}$. Since $\tr(yy^*W)=\langle Wy,y\rangle$, the left-hand side equals $\min\{\langle Wy,y\rangle: y\in S_{\lambda_{\max}},\|y\|=1\}$. In the same way, the
	right-hand side equals  $\max\{\langle Wz,z\rangle: z\in S_{\lambda_{\min}},\|z\|=1\}$,
	contradicting (4).\\
	$(3)\Rightarrow(2)$ Let $Y\in\calD_{S_{\lambda_{\max}}}$ and
	$Z\in\calD_{S_{\lambda_{\min}}}$ such that
	$\Phi(Y)=\Phi(Z)$ and put
	$X=\frac{Y-Z}{2}$. Since $S_{\lambda_{\max}}\perp S_{\lambda_{\min}}$ we have $X_{+}=\frac 12Y$,
	$X_{-}=\frac 12 Z$ and $\|X\|_1=1$, so $X\neq 0$. Moreover,
	\[X=E_+XE_++E_-XE_-,\ \text{with } E_+XE_+=\frac12 Y\geq 0,\ 
	E_-XE_-=-\frac 12Z\leq 0,
	\] 
	and $X\in\calW_1^{\circ}$, since $\Phi(X)=\overline{0}$. 
	Finally $AE_+=\lambda_{\max} E_+$ and
	$AE_-=\lambda_{\min} E_-$ give
	$\tr(AX)=\frac{\lambda_{\max}-\lambda_{\min}}2=\|A\|=\|A\|\,\|X\|_1$.\\	
	$(2)\Rightarrow(1)$
	This implication is immediate, and is \cite[Proposition 5.1]{rieffel}: if there exists $X$ satisfying \eqref{X witness}, then the functional $\varphi(\cdot)=\tr(X\cdot)$ satisfies
	\[\|A\|\,\|X\|_1=\varphi(A)=\tr(AX)=\tr\big((A+W)X\big)\leq\|A+W\|\,\|X\|_1,\]
	for every $W\in\calW$, since $X\in\calW_1^{\circ}$ and $|\tr(RX)|\leq\|R\|\,\|X\|_1$.
\end{proof}
\begin{remark}[Comments on Theorem \ref{teo mini equiv general}]
	\begin{enumerate}
		\item The hypothesis $I\in\calW$ intervenes at exactly two points of the argument. One is through the balanced spectrum
		property, which for the matrices $A$ of the statement holds anyway. The other is the implication $(2)\Rightarrow(3)$, where $\tr(X)=0$ is what allows the two parts of the witness to be normalized to density matrices. Without it one still
		obtains two positive matrices supported in $S_{\lambda_{\max}}$ and in $S_{\lambda_{\min}}$ with the same moments, but not of the same trace. The separation argument and the implications $(3)\Rightarrow(2)\Rightarrow(1)$ do not use it.
		\item $I\in \calW$ implies that every $X\in \calW_1^{\circ}$ is a trace zero matrix. 
		\item Note that condition $X=E_+XE_++E_-XE_-$ is equivalent to state that $X$ equals its pinching operator.
		\item 	When $\calW$ is a $C^{*}$-subalgebra of $M_n(\C)$ containing $I$, the
		equivalence between conditions (1) and (3) follows from \cite[Theorem 2 and Corollary 1]{BV-2026}. The proofs given there do not use the product of $\calW$: it intervenes only
		through the annihilator $\calW_1^{\circ}$ and through $I\in\calW$.	
	\end{enumerate}
\end{remark}
We now consider the general case of a subspace $\calW$ of Hermitian matrices, where possibly $I\notin \calW$.

\begin{notation}
	Let $\calW\subset M_n^h(\C)$ be a subspace, $A_0\in M_n^h(\C)\setminus \calW$ and $W_0\in \calW$. We abbreviate $S_{\mathrm{MAX}}:=S_{\mathrm{MAX}}(A_0+W_0)$ for the subspace introduced in \eqref{SMAX generico}, which by item 1 of Proposition \ref{lema prop smax} can be written as
		\begin{equation}\label{def Smax}
		S_{\mathrm{MAX}}=\operatorname{span}\{v\in \C^n:\ (A_0+W_0)v=\pm\|A_0+W_0\|v\}.
	\end{equation}
\end{notation}
Following \cite{singer}, we define the extremal functionals on $M_n^h(\C)$.
\begin{definition}
	Let $S_1^*:=\left\{f\in \left(M_n^h(\C)\right)^*:\,\|f\|=1\right\}$ denote the unit sphere of the dual space.
	A functional $f\in S_1^*$ is \emph{extremal} if, whenever $f=\frac{f_1+f_2}{2}$ with $f_1,f_2\in S_1^*$, then $f=f_1=f_2$. We denote by $\calE^h\subseteq S_1^*$ the set of extremals. 
\end{definition}
Since every $f\in S_1^*$ can be written as $f(\cdot)=\tr(F\cdot)$ with a unique $F\in  M_n^h(\C)$ with $\|F\|_1=1$, we can identify the set $\calE^h$ with the set of $F\in  M_n^h(\C)$, $\|F\|_1=1$ such that, if there exist Hermitian matrices $F_1,F_2$ with $\|F_1\|_1=\|F_2\|_1=1$ and $F=\frac{F_1+F_2}{2}$, then $F=F_1=F_2$. 
\begin{lemma}\label{lem:extremales}
	Let $u,v \in \mathbb{C}^n$ be unit vectors and let $F = uv^* \in M_n(\C)$.
	If $F$ is self-adjoint ($F = F^*$), then $u = v$ or $u = -v$.
\end{lemma}
\begin{proof}
	Applying the identity $uv^*=vu^*$ to $u$ gives $\langle u,v\rangle\, u=v$, where $\langle x,y\rangle=y^*x$. Hence $v=\lambda u$ with $\lambda=\langle u,v\rangle$, and $|\lambda|=1$ because $\|u\|=\|v\|=1$. Then $\overline{\lambda}\,uu^*=uv^*=vu^*=\lambda\,uu^*$, so that $\lambda\in\R$, that is, $\lambda=\pm1$.
\end{proof}

This result is probably known, but since we could not find a reference, we include a proof.
\begin{proposition}\label{prop:extremales}
	The set of extremals in the unit sphere of $(M_n^h(\C))^*$ is 
	\[\calE^h=\{\pm vv^*: v\in \C^n, \|v\|=1\}.\]
\end{proposition}
\begin{proof}
	Via $F\mapsto \tr(\,\cdot\,F)$, the dual $(M_n^h(\C))^*$ is isometric to $M_n^h(\C)$ with the trace norm $\|\cdot\|_1$. If $\|F\|_1=1$ and $F=\frac{F_1+F_2}{2}$ with $\|F_1\|_1,\|F_2\|_1\leq1$, then $1\leq\frac{\|F_1\|_1+\|F_2\|_1}{2}$ forces $\|F_1\|_1=\|F_2\|_1=1$; hence the extremals are exactly the extreme points of the unit ball of $(M_n^h(\C),\|\cdot\|_1)$, and we may work with convex combinations $tG+(1-t)H$, $0<t<1$, of elements of this ball.
	\\
	Let $F\in\calE^h$, so $\|F\|_1=1$. Writing
	$F=\sum_{\lambda_i\neq0}|\lambda_i|\,(\operatorname{sgn}\lambda_i\, v_iv_i^*)$
	(spectral decomposition, with $\{v_i\}$ orthonormal) exhibits $F$ as a convex combination of the points
	$\pm v_iv_i^*$ of the unit ball. Since $F$ is extreme, it must coincide
	with each of them; therefore $F$ has a single nonzero eigenvalue $\pm1$, and $F=\pm vv^*$ with $\|v\|=1$.
	\\
	Conversely, recall that $\tr G\le\|G\|_1$ for every $G\in M_n^h(\C)$, with equality if and only if $G\ge 0$. Let $\|v\|=1$ and suppose $vv^*=tG+(1-t)H$ with $0<t<1$ and
	$\|G\|_1,\|H\|_1\le1$. Then
	$1=t\tr G+(1-t)\tr H\le t\|G\|_1+(1-t)\|H\|_1\le1$,
	so $\tr G=\|G\|_1=1$, and likewise for $H$. Hence $G,H\ge0$. For
	$w\perp v$ we get $t\langle Gw,w\rangle+(1-t)\langle Hw,w\rangle=\langle vv^*w,w\rangle=0$, so $\langle Gw,w\rangle=0$ and, since $G\geq0$,
	$Gw=0$. Thus $G=\mu vv^*$, and $\tr G=1$ gives $G=vv^*$; likewise
	$H=vv^*$. Finally, if $-vv^*=tG+(1-t)H$, then $vv^*=t(-G)+(1-t)(-H)$, so $-G=-H=vv^*$. Hence $\pm vv^*\in\calE^h$.
\end{proof}

By \cite[Theorem 2]{LR81}, the set of extremals in the unit sphere of $\left( M_n(\C)\right)^*$ is
\[\calE=\{F\in M_n(\C):\ F=uv^*,\ \|u\|=\|v\|=1\}.\]
\begin{remark}
	Comparing both descriptions, $\calE^h$ is exactly the set of self-adjoint elements of $\calE$:
	indeed, by Lemma \ref{lem:extremales} a self-adjoint $uv^*$ with $\|u\|=\|v\|=1$ is of the form
	$\pm vv^*$.
\end{remark}

\begin{theorem}\label{teo minimality merged}
 Let $\calW \subseteq M_n^h(\C)$ be a subspace with $\dim_{\R}\calW = r$, and fix
$A_0 \in M_n^h(\C)\setminus\calW$, $W_0 \in \calW$. Let
$S_{\mathrm{MAX}}$ be as in~\eqref{def Smax}. Then the following are equivalent:	
\begin{enumerate}
		\item $\|A_0+W_0\|=\text{dist}(A_0,\calW)$, that is, $A_0+W_0$ is $\calW$-minimal.
		\item There exist $l$ unit vectors $v_1,\dots,v_l\in S_{\mathrm{MAX}}$, with $1\leq l\leq r+1$, which are moreover eigenvectors of $A_0+W_0$ associated with the eigenvalue $\|A_0+W_0\|$ or with the eigenvalue $-\|A_0+W_0\|$, together with vectors $u_i\in S_{\mathrm{MAX}}$ and weights $0\leq \omega_i\leq 1$ with $\sum_{i=1}^{l}\omega_i=1$, such that
		\begin{equation}\label{cond1}
			\sum_{i=1}^{l}\omega_i \langle Wv_i, u_i\rangle=0, \text{ for every } W\in \calW,
		\end{equation}
		where the vectors $u_i$ are given by
		\begin{equation} \label{cond2}
			\begin{array}{lll} 
				u_i=v_i&\text{if}&(A_0+W_0)v_i=\|A_0+W_0\|v_i,\\
				u_i=-v_i&\text{if}&(A_0+W_0)v_i=-\|A_0+W_0\|v_i.
			\end{array}
		\end{equation}
		
		\item 
		\begin{equation} \label{0 in moments}
			\left\lbrace \begin{array}{lcl}
				\overline{0}\in m_{S_{\mathrm{MAX}}}&\text{if }& \lambda_{\max}+\lambda_{\min}\neq 0\\
				\overline{0}\in \operatorname{conv}\left(  m_{S_{\lambda_{\max}}}\cup - m_{S_{\lambda_{\min}}}\right) & \text{if } &\lambda_{\max}+\lambda_{\min}= 0
			\end{array}\right. 
		\end{equation}
		where $\lambda_{\max}=\max\sigma(A_0+W_0)$ and $\lambda_{\min}=\min\sigma(A_0+W_0)$.
		\item $\overline{0}\in \partial(\|A(\cdot)\|)(x)$, where $x=(x_1,x_2,\dots,x_r)\in \R^r$ is such that $A(x)=A_0+\sum_{i=1}^{r}x_iB^i=A_0+W_0$.
		\item There exists $F\in \calW^{\circ}_1$ such that $F\in \partial\|A_0+W_0\| $, i.e.:
		\[\tr(F(A_0+W_0))=\tr(FA_0)=\|A_0+W_0\|.\]
	\end{enumerate}		
\end{theorem}

\begin{proof}
$(1)\Rightarrow(2)$ By Theorem 1.1 in \cite{singer}, page 170, $-W_0$ is a best approximation of $A_0$ in $\calW$ if and only if there exist $l$ extremal functionals $F_i\in \calE^h$ and weights $0\leq \omega_i\leq 1$ with $\sum_{i=1}^{l}\omega_i=1$ such that
	\begin{equation}
		\begin{split}\label{eqsinger}
		\sum_{i=1}^{l}\omega_i\tr(F_iW)&=0, \ \forall W\in \calW\\
		\tr(F_i(A_0+W_0))&=\|A_0+W_0\|, \text{ for every } 1\leq i\leq l.
		\end{split}
	\end{equation}
	By Proposition \ref{prop:extremales}, each $F_i$ is of the form $F_i=\pm v_iv_i^*$ with $\|v_i\|=1$, and $\tr(v_iv_i^*R)=\langle Rv_i,v_i\rangle$ for any matrix $R$; then \eqref{eqsinger} is equivalent to 
	\begin{equation*}
		\begin{split}
			\sum_{i=1}^{l}\omega_i \langle Wv_i, u_i\rangle&=0, \ \forall W\in \calW\\
			\tr(F_i(A_0+W_0))&= u_i^*(A_0+W_0)v_i=\|A_0+W_0\|, \text{ for every } 1\leq i\leq l,
		\end{split}
	\end{equation*}
	where each $v_i\in S_{\mathrm{MAX}}$ and $u_i=\pm v_i$. 
	\begin{itemize}
		\item If $u_i=v_i$, then $v_i$ is an eigenvector of $A_0+W_0$ associated with $\|A_0+W_0\|$.		
		\item If $u_i=-v_i$, then $v_i$ is an eigenvector of $A_0+W_0$ associated with $-\|A_0+W_0\|$.
	\end{itemize}	
	$(2)\Rightarrow(1)$ If there exist $l$ unit vectors $v_i,u_i\in S_{\mathrm{MAX}}$ and weights $0\leq \omega_i\leq 1$ with $\sum_{i=1}^{l}\omega_i=1$ satisfying \eqref{cond1} and \eqref{cond2}, then by \eqref{cond2} the matrix $F=\sum_{i=1}^{l}\omega_i u_iv_i^*$ satisfies
	\[\begin{cases*}
	F\in M_n^h(\C), \\
	\tr(FW)=0, \ \forall\ W\in \calW,\\
	\tr(FA_0)=\tr(F(A_0+W_0))=\|A_0+W_0\|.
	\end{cases*}\]
	Moreover, $\|F\|_1\leq \sum_{i=1}^{l}\omega_i=1$, while from $\tr(F(A_0+W_0))=\|A_0+W_0\|$ we get $\|F\|_1\geq 1$; hence $\|F\|_1=1$, and therefore $F\in\calW^{\circ}_1\cap\partial\|A_0+W_0\|$. Consequently, since $\tr(FW)=0$ for every $W\in\calW$,
	\begin{align*}
	\|A_0+W_0\|&=\tr\left( F(A_0+W_0)\right)=\tr\left( F(A_0+W)\right)\\
	&\leq\|A_0+W\|\,\|F\|_1=\|A_0+W\|,
	\end{align*}
	for every $W\in\calW$; that is, $A_0+W_0$ is $\calW$-minimal.

	The same computation proves $(5)\Rightarrow(1)$, while $(1)\Rightarrow(5)$ follows by combining $(1)\Rightarrow(2)$ with the construction of $F$ above. This gives a self-contained proof of
	\[\|A_0+W_0\|=\min_{W\in \calW}\|A_0+W\| \Leftrightarrow 	\exists F\in \partial\|A_0+W_0\|\cap \calW^{\circ}_1,   
	\]
	an equivalence which is also a particular case of \cite[Corollary 1]{zietak-1996}.
	
	$(2)\Rightarrow(3)$	Since $\{B^j\}_{j=1}^r$ is a basis of $\calW$, condition \eqref{cond1} for any $W\in \calW$ is equivalent to
	\[\sum_{i=1}^{l}\omega_i \langle B^jv_i,u_i\rangle= 0,\ \text{for all } 1\leq j\leq r,\]
	which implies that
	\begin{equation}\label{eq =0}
		\sum_{i=1}^{l}\omega_i \left( \langle B^1v_i,u_i\rangle,\cdots,  \langle B^rv_i,u_i\rangle\right) = (0,\dots,0),
	\end{equation}
	that is, $\overline{0}$ is a convex combination of the vectors $\left( \langle B^1v_i,u_i\rangle,\cdots,  \langle B^rv_i,u_i\rangle\right)$, $1\leq i\leq l$.
	We distinguish three cases:
	\begin{itemize}
		\item If $\|A_0+W_0\|\in \sigma(A_0+W_0)$ and $-\|A_0+W_0\|\notin \sigma(A_0+W_0)$, then  $S_{\mathrm{MAX}}=S_{\lambda_{\max}}$, so $v_i\in S_{\lambda_{\max}}$ and $\langle B^jv_i,u_i\rangle=\langle B^jv_i,v_i\rangle$ for every $1\leq j\leq r$.  Therefore by \eqref{eq =0}, $\overline{0}\in m_{S_{\lambda_{\max}}}= m_{S_{\mathrm{MAX}}}$.
		\item If $-\|A_0+W_0\|\in \sigma(A_0+W_0)$ and $\|A_0+W_0\|\notin \sigma(A_0+W_0)$, then  $S_{\mathrm{MAX}}=S_{\lambda_{\min}}$, so $v_i\in S_{\lambda_{\min}}$ and $\langle B^jv_i,u_i\rangle=-\langle B^jv_i,v_i\rangle$ for every $1\leq j\leq r$. Therefore, by \eqref{eq =0}, $\overline{0}\in -m_{S_{\lambda_{\min}}}= -m_{S_{\mathrm{MAX}}}$ which is equivalent to $\overline{0}\in m_{S_{\mathrm{MAX}}}$.
		\item If $\pm\|A_0+W_0\|\in\sigma(A_0+W_0)$, then $S_{\mathrm{MAX}}=S_{\lambda_{\min}}\oplus S_{\lambda_{\max}}$. By hypothesis, each $v_i$ is an eigenvector of $A_0+W_0$ corresponding to $\|A_0+W_0\|$ or to $-\|A_0+W_0\|$; after reindexing, we may assume that
		\[\begin{array}{cccl}
			(A_0+W_0)v_i=\|A_0+W_0\|v_i,&v_i\in S_{\lambda_{\max}}& \text{ for every } 1\leq i\leq d,\\
			(A_0+W_0)v_i=-\|A_0+W_0\|v_i&v_i\in S_{\lambda_{\min}}& \text{ for every } d+1\leq i\leq l,\\
		\end{array}\]
		 and, writing $x_i=\left( \langle B^1v_i,v_i\rangle,\dots,\langle B^rv_i,v_i\rangle\right)$ for $1\leq i\leq l$, condition \eqref{cond1} gives
		\[
		\sum_{i=1}^{d}\omega_i x_i+ \sum_{i=d+1}^{l}\omega_i \left( \langle B^jv_i,-v_i\rangle\right)_{j=1}^{r}=\sum_{i=1}^{d}\omega_i x_i-\sum_{i=d+1}^{l}\omega_i x_i=\overline{0}.
		\]
		Then,
		\[\overline{0}\in \operatorname{conv}\left( m_{S_{\lambda_{\max}}}\cup -m_{S_{\lambda_{\min}}} \right),\]
		since $x_i\in m_{S_{\lambda_{\max}}}$ for $1\leq i\leq d$ and $-x_i\in -m_{S_{\lambda_{\min}}}$ for $d+1\leq i\leq l$.
		
	\end{itemize}
	$(3)\Rightarrow(2)$ First suppose that $\lambda_{\max}+\lambda_{\min}\neq 0$, that is $S_{\mathrm{MAX}}=S_{\lambda_{\max}}$ or $S_{\mathrm{MAX}}=S_{\lambda_{\min}}$. Since $\overline{0}\in m_{S_{\mathrm{MAX}}}$ and $m_{S_{\mathrm{MAX}}}$ is the convex hull of $\{(\langle B^1v,v\rangle,\cdots,\langle B^rv,v\rangle): v\in S_{\mathrm{MAX}},\ \|v\|=1\}\subseteq\R^r$, using Carath\'eodory's theorem (see for instance \cite[Theorem 17.1]{Rockafellar}) there exist $l\leq r+1$ unit vectors $v_i\in S_{\mathrm{MAX}}$ and weights $0\leq \alpha_i\leq 1$ with $\sum_{i=1}^{l}\alpha_{i}=1$ such that  $y_i=(\langle B^1v_i,v_i\rangle,\cdots, \langle B^rv_i,v_i\rangle)$ for every $i$ and
	\begin{align*}
		\overline{0}&=\sum_{i=1}^{l}\alpha_iy_i=\sum_{i=1}^{l}\alpha_i (\langle B^1v_i,v_i\rangle,\cdots, \langle B^rv_i,v_i\rangle).
	\end{align*}
	This implies that $0= \sum_{i=1}^{l}\alpha_i \langle B^jv_i,v_i\rangle$, for every $1\leq j\leq r$. Therefore, every $W\in \calW$ satisfies
	\[0= \sum_{i=1}^{l}\alpha_i \langle Wv_i,v_i\rangle.\]
	We now distinguish two cases.
	\begin{itemize}
		\item If $S_{\mathrm{MAX}}=S_{\lambda_{\max}}$, then $(A_0+W_0)v_i=\|A_0+W_0\|v_i$ for every $1\leq i\leq l$.	
		\item If $S_{\mathrm{MAX}}=S_{\lambda_{\min}}$, then $(A_0+W_0)v_i=-\|A_0+W_0\|v_i$ for every $1\leq i\leq l$, so that $u_i=-v_i$ by \eqref{cond2}. Hence, since $0=\sum_{i=1}^{l}\alpha_i\langle Wv_i,v_i\rangle$,
		\[0= \sum_{i=1}^{l}\alpha_i \langle Wv_i,-v_i\rangle.\]
	\end{itemize}
	In both cases, conditions \eqref{cond1} and \eqref{cond2} are satisfied.

	Now suppose that $\lambda_{\max}+\lambda_{\min}=0$, that is $S_{\mathrm{MAX}}=S_{\lambda_{\max}}\oplus S_{\lambda_{\min}}$. Set
	\begin{align*}
		T_{\max}&=\{(\langle B^1v,v\rangle,\cdots,\langle B^rv,v\rangle): v\in S_{\lambda_{\max}},\ \|v\|=1\}\quad\text{and}\\
		T_{\min}&=\{- (\langle B^1v,v\rangle,\cdots,\langle B^rv,v\rangle): v\in S_{\lambda_{\min}},\ \|v\|=1\}.
	\end{align*}
	Since $m_{S_{\lambda_{\max}}}=\operatorname{conv}(T_{\max})$, $-m_{S_{\lambda_{\min}}}=\operatorname{conv}(T_{\min})$ and $\operatorname{conv}(C\cup D)=\operatorname{conv}(\operatorname{conv}(C)\cup\operatorname{conv}(D))$
	for any $C,D\subseteq\R^r$, the condition
	\[\overline{0}\in \operatorname{conv}\left( m_{S_{\lambda_{\max}}}\cup -m_{S_{\lambda_{\min}}} \right)\quad\text{ is equivalent to } \overline{0}\in\operatorname{conv}(T_{\max}\cup T_{\min}).\] 
	By Carath\'eodory's theorem in $\R^r$, $\overline{0}$ is a convex combination of at most $r+1$ points of $T_{\max}\cup T_{\min}$, that is,
	\[\sum_{i=1}^{l}\alpha_{i}y_i=\overline{0}\]
	for some $l\leq r+1$, weights $\alpha_i\geq 0$ with $\sum_{i=1}^{l}\alpha_i=1$, and $y_i\in T_{\max}\cup T_{\min}$. Each $y_i$ is of the form
	\[y_i=s_i(\langle B^1w_i,w_i\rangle,\cdots,\langle B^rw_i,w_i\rangle),\]
	where $w_i$ is a unit vector, $s_i=1$ if $w_i\in S_{\lambda_{\max}}$ and $s_i=-1$ if $w_i\in S_{\lambda_{\min}}$. Therefore,
	\[\sum_{i=1}^{l}\alpha_{i}s_i\langle B^jw_i,w_i\rangle =0, \text{ for every } 1\leq j\leq r.\]
	  Then, regrouping the points of each type and normalizing the weights on each side, there exist $0\leq \alpha\leq 1$, $\rho\in \calD_{S_{\lambda_{\max}}}$ of rank at most $d$ and $\sigma\in \calD_{S_{\lambda_{\min}}}$ of rank at most $k$, with $d+k\leq r+1$, such that $\Phi(\rho)\in m_{S_{\lambda_{\max}}}$, $\Phi(\sigma)\in m_{S_{\lambda_{\min}}}$ (see \eqref{defi phi}) and
	\[\overline{0}=\alpha \Phi(\rho)+(1-\alpha)\Phi(-\sigma)=\Phi\left(\alpha \rho-(1-\alpha)\sigma \right) .\]
	Thus, $\tr\left(\left(\alpha \rho-(1-\alpha)\sigma \right)B^j \right)=0 $ for every $1\leq j\leq r$, and this implies that  
	\[\tr\left(\left(\alpha \rho-(1-\alpha)\sigma \right)W \right)=0\]
	for all $W\in \calW$. Since $\rho$ and $\sigma$ have orthogonal supports, $\|\alpha \rho-(1-\alpha)\sigma\|_1=\alpha+(1-\alpha)=1$, and therefore $\alpha \rho-(1-\alpha)\sigma\in \calW_1^{\circ}$.
	
	On the other hand, since $\rho\in \calD_{S_{\lambda_{\max}}}$ and $\sigma\in \calD_{S_{\lambda_{\min}}}$, both can be written as
	\[\rho=\sum_{i=1}^{d}\beta_i v_iv_i^* \ \text{ and } \sigma=\sum_{j=1}^{k}\gamma_j z_jz_j^*,\]
	with $0\leq \beta_i,\gamma_j\leq 1$,  $\sum_{i=1}^{d}\beta_i=\sum_{j=1}^{k}\gamma_j=1$, and orthonormal sets $\{v_i\}_{i=1}^d\subset S_{\lambda_{\max}}$, $\{z_j\}_{j=1}^k\subset S_{\lambda_{\min}}$. Then, 
	\begin{align*}
		\tr\left(\left(\alpha \rho-(1-\alpha)\sigma \right)W \right)&=\alpha\tr(\rho W)-(1-\alpha)\tr(\sigma W) \\
		&=\alpha\tr\left(\sum_{i=1}^{d}\beta_i v_iv_i^*W\right)-(1-\alpha)\tr\left(\sum_{j=1}^{k}\gamma_j z_jz_j^*W\right) \\
		&=\alpha\sum_{i=1}^{d}\beta_i\tr\left( v_iv_i^*W\right)-(1-\alpha)\sum_{j=1}^{k}\gamma_j\tr\left( z_jz_j^*W\right) \\
		&=\alpha\sum_{i=1}^{d}\beta_i\langle Wv_i,v_i\rangle-(1-\alpha)\sum_{j=1}^{k}\gamma_j\langle Wz_j,z_j\rangle \\
		&=\alpha\sum_{i=1}^{d}\beta_i\langle Wv_i,v_i\rangle+(1-\alpha)\sum_{j=1}^{k}\gamma_j\langle Wz_j,-z_j\rangle \\		
		&=\sum_{i=1}^{d}\alpha\beta_i\langle Wv_i,v_i\rangle+\sum_{j=1}^{k}(1-\alpha)\gamma_j\langle Wz_j,-z_j\rangle.		
	\end{align*}
	Observe that $\alpha\sum_{i=1}^{d}\beta_i+(1-\alpha)\sum_{j=1}^{k}\gamma_j= \alpha+1-\alpha=1$. Then if we define
	\[\omega_i=\left\lbrace \begin{array}{lll}
		\alpha\beta_i&\text{if}&1\leq i\leq d\\
		(1-\alpha)\gamma_{i-d}&\text{if}&d+1\leq i\leq d+k,\\
	\end{array}\right. \]
	we obtain that $0\leq \omega_i\leq 1$ for every $i$ and $\sum_{i=1}^{d+k}\omega_i=1$. Now writing $v_{d+j}:=z_j$ for $1\leq j\leq k$, and $u_i:=v_i$ (resp.\ $u_i:=-v_i$) for $1\leq i\leq d$ (resp.\ $d+1\leq i\leq d+k$), we get
	\[\sum_{i=1}^{d+k}\omega_i\langle Wv_i,u_i\rangle=0,\]
	which are exactly conditions \eqref{cond1} and \eqref{cond2}, with $l=d+k\leq r+1$ by the Carath\'eodory decomposition above. Moreover $v_i\in S_{\lambda_{\max}}$ for $1\leq i\leq d$ and $v_i\in S_{\lambda_{\min}}$ for $d+1\leq i\leq d+k$, so that all of them are unit vectors of $S_{\mathrm{MAX}}$ and eigenvectors of $A_0+W_0$ associated with $\|A_0+W_0\|$ or with $-\|A_0+W_0\|$, as required in item (2).
	This completes the proof of $(3)\Rightarrow(2)$.

	$(3)\Leftrightarrow(4)$ By Corollary \ref{formula for subdif composition norm}, each condition in \eqref{0 in moments} is equivalent to $\overline{0}\in\partial(\|A(\cdot)\|)(x)$, since
		\begin{itemize}
		\item If $\|A(x)\|=\lambda_{\max}(A(x))> -\lambda_{\min}(A(x))$, then $$\partial(\|A(\cdot)\|)(x)=m_{S_{\mathrm{MAX}}(A(x))}=m_{S_{\mathrm{MAX}}}.$$
		\item If $\|A(x)\|=-\lambda_{\min}(A(x))> \lambda_{\max}(A(x))$, then $$\partial(\|A(\cdot)\|)(x)=-m_{S_{\mathrm{MAX}}(A(x))}=-m_{S_{\mathrm{MAX}}}.$$
		\end{itemize}
		In both cases, we deduce that $\overline{0}\in \partial(\|A(\cdot)\|)(x)\Leftrightarrow \overline{0}\in m_{S_{\mathrm{MAX}}}$.
		\begin{itemize}
		\item If $\|A(x)\|=-\lambda_{\min}(A(x))=\lambda_{\max}(A(x))$, then 
		$$\partial(\|A(\cdot)\|)(x)={\rm conv}\left(m_{S_{\lambda_{\max}}}\cup-m_{S_{\lambda_{\min}}} \right).$$ 
		Therefore, $\overline{0}\in \partial(\|A(\cdot)\|)(x)\Leftrightarrow \overline{0}\in {\rm conv}\left(m_{S_{\lambda_{\max}}}\cup-m_{S_{\lambda_{\min}}} \right)$.
	\end{itemize}
	
\end{proof}

In the balanced case, condition \eqref{0 in moments} always follows from the intersection condition of Theorem \ref{teo mini equiv general}: if $\Phi(\rho)=\Phi(\sigma)$ for some $\rho\in\calD_{S_{\lambda_{\max}}}$ and $\sigma\in\calD_{S_{\lambda_{\min}}}$, then $\overline{0}=\frac12\Phi(\rho)-\frac12\Phi(\sigma)\in\operatorname{conv}\left(m_{S_{\lambda_{\max}}}\cup -m_{S_{\lambda_{\min}}}\right)$. Both conditions are equivalent when $I\in\calW$, but the converse implication may fail when $I\notin\calW$, as the following remark shows.

\begin{remark}\label{rem:sin identidad}
	When $\lambda_{\max}+\lambda_{\min}=0$, the second condition in \eqref{0 in moments} means that there exist $\alpha\in[0,1]$, $\rho\in\calD_{S_{\lambda_{\max}}}$ and $\sigma\in\calD_{S_{\lambda_{\min}}}$ such that $\alpha\Phi(\rho)-(1-\alpha)\Phi(\sigma)=\overline{0}$. Equivalently,
	$\alpha\,\Phi(\rho)=(1-\alpha)\,\Phi(\sigma)$, that is,
	\[
	\alpha\, m_{S_{\lambda_{\max}}}\cap(1-\alpha)\, m_{S_{\lambda_{\min}}}\neq\emptyset .
	\]	
	By linearity $\Phi(\alpha\rho-(1-\alpha)\sigma)=\overline{0}$, so $F=\alpha\rho-(1-\alpha)\sigma\in\ker(\Phi)=\calW^{\circ}$; moreover, since $\rho$ and $\sigma$ have orthogonal supports, $\|F\|_1=\alpha+(1-\alpha)=1$, so $F\in\calW_1^{\circ}$ (it is the witness of item (5) of Theorem \ref{teo minimality merged}). Therefore
	\[\tr(FW)=0\]
	for every $W\in \calW$. In particular, if $I\in\calW$,
	\[\tr(F)=\alpha\tr(\rho)-(1-\alpha)\tr(\sigma)=2\alpha-1=0.\]
	Hence $\alpha=\frac12$ and $\Phi(\rho)=\Phi(\sigma)\in m_{S_{\lambda_{\max}}}\cap m_{S_{\lambda_{\min}}}$, which is item (3) of Theorem \ref{teo mini equiv general}. If $I\notin\calW$, the weight $\alpha$ is not determined in this way, and a $\calW$-minimal matrix with balanced spectrum may satisfy
	\[
	m_{S_{\lambda_{\max}}}\cap m_{S_{\lambda_{\min}}}=\emptyset
	\quad\text{and}\quad
	m_{S_{\lambda_{\max}}}\cap -m_{S_{\lambda_{\min}}}=\emptyset .
	\]
	For instance, let $n=2$, $\calW=\R B$ with $B=\frac{1}{\sqrt5}\begin{psmallmatrix}1&0\\0&2\end{psmallmatrix}$ (so that $I\notin\calW$), $A_0=\begin{psmallmatrix}1&0\\0&-1\end{psmallmatrix}$ and $W_0=0$. For $t\in\R$,
	\[
	\Big\|A_0+t\sqrt5\, B\Big\|=\max\{|1+t|,\ |2t-1|\}\geq 1,
	\]
	since the first term is greater than $1$ if $t>0$ and the second one is greater than $1$ if $t<0$. Hence $A_0$ is $\calW$-minimal, with $\lambda_{\max}=1=-\lambda_{\min}$. Both eigenvalues are simple, $S_{\lambda_{\max}}=\operatorname{span}\{e_1\}$ and $S_{\lambda_{\min}}=\operatorname{span}\{e_2\}$, so
	\[
	m_{S_{\lambda_{\max}}}=\Big\{\tfrac{1}{\sqrt5}\Big\},\qquad m_{S_{\lambda_{\min}}}=\Big\{\tfrac{2}{\sqrt5}\Big\},
	\]
	and both intersections above are empty. However, $\overline0=\frac23\cdot\frac{1}{\sqrt5}+\frac13\cdot\big(-\frac{2}{\sqrt5}\big)\in\operatorname{conv}\big(m_{S_{\lambda_{\max}}}\cup -m_{S_{\lambda_{\min}}}\big)$, as predicted by Theorem \ref{teo minimality merged}. Here the witness is $F=\frac23e_1e_1^*-\frac13e_2e_2^*$.
	The same phenomenon occurs in Example \ref{ex:2x2 block}, where $\calW$ comes from a completion pattern: there $\alpha=\frac{5+2\sqrt5}{10}\neq\frac12$, and the moments are the single points $(a,0,0,0,0)$ and $(b,0,0,0,0)$ with $a=\frac{5-2\sqrt5}{10}$ and $b=\frac{5+2\sqrt5}{10}$, which are distinct and positive.
\end{remark}

\section{Examples and applications} \label{sec:examples}

\begin{example}\label{ex:2x2 block}
	Let $A_0=\begin{pmatrix}
		0&1&0\\
		1&2&0\\
		0&0&1
	\end{pmatrix}$. Then
	\begin{gather*}
		I_0=\{(1,1),(1,3),(2,3),(3,1),(3,2)\},\\
		\calW=\operatorname{span}\{B^{1,1},B^{1,3},B^{2,3}, B^{3, 1}, B^{3, 2}\}.
	\end{gather*}
	By Proposition \ref{A0 Y comparten ceros}, and since the trace norm of a block diagonal matrix is the sum of the trace norms of its blocks,
	\[\calW^{\circ}_1 =\left\{Y=\begin{pmatrix}
		0&Y_{12}&0\\
		\overline{Y_{12}}&Y_{22}&0\\
		0&0&Y_{33}
	\end{pmatrix}:\ \sqrt{4|Y_{12}|^2+Y_{22}^2}+|Y_{33}|=1\right\},\]
	because the eigenvalues $\frac{1}{2}\left( Y_{22}\pm\sqrt{Y_{22}^2+4|Y_{12}|^2}\right)$ of the $2\times 2$ block have opposite signs, so that its trace norm is their difference. By the duality formula \eqref{BDF}, $\min_{W\in \calW}\|A_0+W\|$ equals
	\begin{equation}\label{opt problem2}
		\max_{(x,y,z)\in \C\times \R\times \R}|2\operatorname{Re}(x)+2y+z|, \quad \text{such that } \sqrt{4|x|^2+y^2}+|z|=1,
	\end{equation}
	where $Y$ has entries $Y_{12}=x$, $Y_{22}=y$, $Y_{33}=z$. The objective does not depend on $\operatorname{Im}(x)$, which only increases $|x|$ in the restriction, so we may assume $x\in\R$; replacing $F$ by $-F$ if necessary, also $2x+2y+z\geq0$, and then $z\geq0$ at a maximizer, since the restriction only involves $|z|$. Put $s=\sqrt{4x^2+y^2}=1-z$. By the Cauchy--Schwarz inequality $2x+2y\leq\sqrt5\, s$, hence $2x+2y+z\leq 1+(\sqrt5-1)s\leq \sqrt5$, with equality if and only if $z=0$ and $(2x,y)=\frac{1}{\sqrt5}(1,2)$. Therefore
	\begin{equation}\label{sol opt problem2}
		\text{dist}(A_0,\calW)=\sqrt{5},
		\qquad \text{attained at}\quad
		x=\frac{1}{2\sqrt{5}},\quad
		y=\frac{2}{\sqrt{5}},\quad
		z=0.
	\end{equation}
	The minimum is attained at $W_0=-2B^{1,1}$, that is, $2B^{1,1}$ is a best approximation of $A_0$ in $\calW$: the eigenvalues of $A_0+W_0$ are $\{1,\sqrt{5},-\sqrt{5}\}$.

	Let us check this with item (3) of Theorem \ref{teo minimality merged}. Since $\lambda_{\max}+\lambda_{\min}=0$, the condition to be verified is the second one in \eqref{0 in moments}. Both eigenvalues are simple, with $S_{\lambda_{\max}}=\operatorname{span}\{(\sqrt{5}-2,1,0)\}$ and $S_{\lambda_{\min}}=\operatorname{span}\{(-(2+\sqrt{5}),1,0)\}$, so that $\calD_{S_{\lambda_{\max}}}=\{\rho_+\}$ and $\calD_{S_{\lambda_{\min}}}=\{\rho_-\}$, where
	\[\rho_+=\frac{1}{10}\begin{pmatrix}
		5-2\sqrt{5} & \sqrt{5} & 0\\
		\sqrt{5} & 5+2\sqrt{5} & 0\\
		0 & 0 & 0
	\end{pmatrix},
	\qquad
	\rho_-=\frac{1}{10}\begin{pmatrix}
		5+2\sqrt{5} & -\sqrt{5} & 0\\
		-\sqrt{5} & 5-2\sqrt{5} & 0\\
		0 & 0 & 0
	\end{pmatrix}.\]
	Since $\rho_\pm$ vanish at the positions $(1,3)$ and $(2,3)$, by \eqref{Yij} and \eqref{Yii} only the coordinate corresponding to $B^{1,1}$ survives: $m_{S_{\lambda_{\max}}}=\{(a,0,0,0,0)\}$ and $-m_{S_{\lambda_{\min}}}=\{(-b,0,0,0,0)\}$, with $a=(\rho_+)_{11}$ and $b=(\rho_-)_{11}$. As $\rho_++\rho_-$ is the orthogonal projection onto $\operatorname{span}\{e_1,e_2\}$, we have $a+b=1$, so $\alpha a-(1-\alpha)b=0$ for
	\begin{equation}\label{peso ejemplo 2}
		\alpha=b=\frac{5+2\sqrt5}{10}\in(0,1).
	\end{equation}
	Hence $\overline{0}\in \operatorname{conv}\left( m_{S_{\lambda_{\max}}}\cup -m_{S_{\lambda_{\min}}}\right)$, as expected. Note that, since $a\neq b$, $m_{S_{\lambda_{\max}}}\cap m_{S_{\lambda_{\min}}}=\emptyset$: here $I\notin\calW$ (see Remark \ref{rem:sin identidad}). Moreover, the same weight recovers the maximizer of \eqref{opt problem2}: $F=\alpha\rho_+-(1-\alpha)\rho_-$ has $F_{11}=F_{33}=0$, $F_{12}=\frac{1}{2\sqrt5}$ and $F_{22}=\frac{2}{\sqrt5}$, which is \eqref{sol opt problem2}. This illustrates the equivalence between items (3) and (5) of Theorem \ref{teo minimality merged}.
\end{example}

\begin{figure}[ht]
	\centering
	\begin{tikzpicture}[scale=3]
		\fill[blue!15] (0.5,0) ellipse [x radius=0.5, y radius=0.7071];
		\draw[-latex, gray] (-0.1,0) -- (1.3,0) node[right, black] {$B^{1,1}$};
		\draw[-latex, gray] (0,-0.85) -- (0,0.85) node[above, black] {$B^{1,2}$};
		\draw[thick, blue!50!black] (0.5,0) ellipse [x radius=0.5, y radius=0.7071];
		\draw[gray] (1,0) -- (1,-0.03) node[below right, black, inner sep=1pt] {\scriptsize $1$};
		\fill[red!75!black] (0,0) circle (0.6pt) node[below left, black] {\small $\overline 0$};
		\node[blue!50!black] at (0.55,0.22) {$m_{S_{\mathrm{MAX}}}$};
	\end{tikzpicture}
	\caption{Example \ref{ex:no unicidad} with $\alpha=2$. Projection of $m_{S_{\mathrm{MAX}}}$ onto the coordinates corresponding to $B^{1,1},B^{1,2}$, which is $\{(\beta,\sqrt2\operatorname{Re}c):\ (\beta-\frac12)^2+(\operatorname{Re}c)^2\leq\frac14\}$. The origin $\overline 0=\Phi(e_2e_2^*)$ belongs to the moment, but only as a boundary point.}
	\label{fig:momento elipse}
\end{figure}
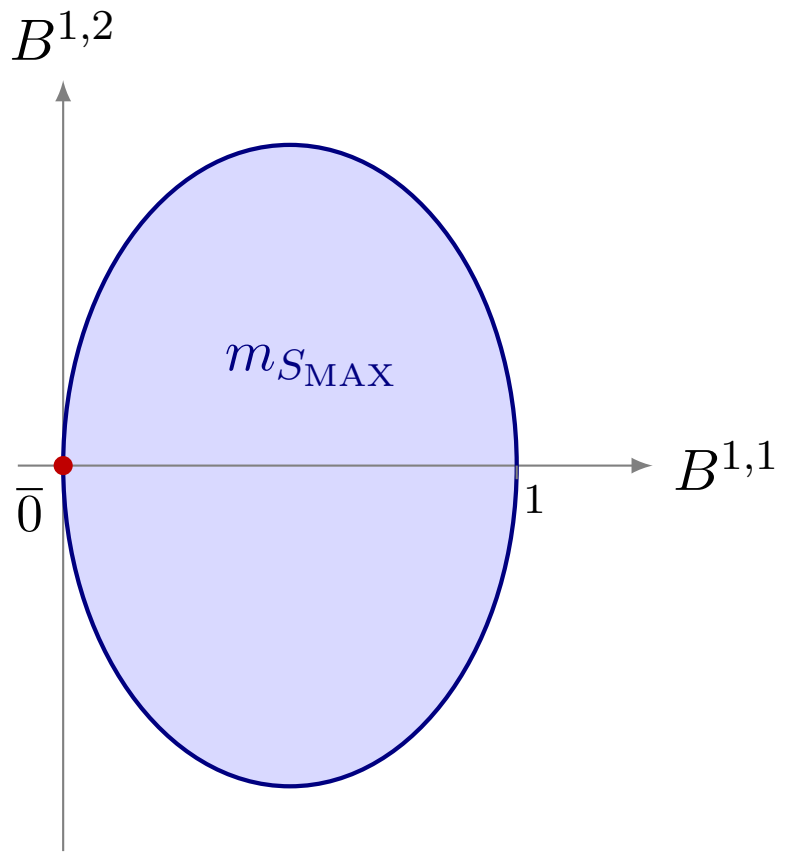

\begin{example}\label{ex:no unicidad}
	This example shows that best approximations need not be unique, and how $S_{\mathrm{MAX}}$ and its moment change along the set of best approximations. Let $A_0=\Diag(0,2,1)$. Then
	\begin{gather*}
		I_0=\{(1,1),(1,2),(1,3),(2,1),(2,3),(3,1),(3,2)\},\\
		\calW=\operatorname{span}\{B^{1,1},B^{1,2}, B^{1,3},B^{2,1},B^{2,3}, B^{3, 1}, B^{3, 2}\},
	\end{gather*}
	so that the only specified entries are the diagonal entries $(2,2)$ and $(3,3)$. For every $W\in\calW$,
	\[\|A_0+W\|\geq|(A_0+W)_{22}|=2=\|A_0\|,\]
	hence $\text{dist}(A_0,\calW)=2$, $A_0$ is $\calW$-minimal and $0$ is a best approximation of $A_0$ in $\calW$. It is not unique: for every $-2\leq\alpha\leq2$, $A_0+\alpha B^{1,1}=\Diag(\alpha,2,1)$ has norm $2$, so $-\alpha B^{1,1}$ is also a best approximation. Along this segment the different situations of Theorem \ref{teo minimality merged} occur:
	\begin{itemize}
		\item If $|\alpha|<2$, then $\lambda_{\max}=2>-\lambda_{\min}$, $S_{\mathrm{MAX}}=\operatorname{span}\{e_2\}$ and $m_{S_{\mathrm{MAX}}}=\{\overline 0\}$, because $\langle B^ke_2,e_2\rangle=(B^k)_{22}=0$ for every $k$.
		\item If $\alpha=2$, then still $\lambda_{\max}=2>-\lambda_{\min}$, but now $S_{\mathrm{MAX}}=\operatorname{span}\{e_1,e_2\}$ and every $\rho\in\calD_{S_{\mathrm{MAX}}}$ is of the form
		\[
		\rho=\begin{pmatrix}\beta&c\\ \bar c&1-\beta\end{pmatrix}\oplus 0,
		\qquad 0\leq\beta\leq1,\ \ |c|^2\leq\beta(1-\beta).
		\]
		Since $\tr(\rho B^{1,1})=\beta$, $\tr(\rho B^{1,2})=\sqrt2\,\operatorname{Re}c$, $\tr(\rho B^{2,1})=-\sqrt2\,\operatorname{Im}c$,
		and $\tr(\rho B^k)=0$ for the remaining elements of $\B$ (which involve the third row or column),
		item \ref{mS = diag DsubS} of Proposition \ref{prop momento} gives, in the coordinates corresponding to
		$B^{1,1},B^{1,2},B^{2,1}$ (all other coordinates being zero),
		\[
		m_{S_{\mathrm{MAX}}}=\bigl\{(\beta,\sqrt2\,\operatorname{Re}c,-\sqrt2\,\operatorname{Im}c):\ 0\leq\beta\leq1,\ |c|^2\leq\beta(1-\beta)\bigr\},
		\]
		a solid ellipsoid. It contains $\overline 0=\Phi(e_2e_2^*)$, as required by condition (3) of Theorem \ref{teo minimality merged}, but only as a boundary point (see Figure \ref{fig:momento elipse}).
		\item If $\alpha=-2$, the spectrum is balanced: $S_{\lambda_{\max}}=\operatorname{span}\{e_2\}$, $S_{\lambda_{\min}}=\operatorname{span}\{e_1\}$, and the second condition of \eqref{0 in moments} holds since $m_{S_{\lambda_{\max}}}=\{\overline 0\}$.
	\end{itemize}
\end{example}

\subsection{Adjacency matrices of graphs}\label{sec:graphs}

Let $G$ be a connected simple graph with vertices $\{1,\dots,n\}$, $n\geq2$, and let $A_0=A(G)$ be its adjacency matrix. Then $I_0$ consists of the non-edges of $G$ together with the diagonal, and by Proposition \ref{A0 Y comparten ceros}, $\calW^\circ_1$ consists of the Hermitian matrices of trace norm one supported on the edges of $G$; in particular, they have zero diagonal. Completing $A_0$ within $\calW$ means choosing freely the weights of the non-edges and of the diagonal (that is, adding loops), and the natural question is when $A(G)$ itself is a best completion. We shall use the following facts.

\begin{lemma}[Perron--Frobenius for connected graphs]\label{lem:PF}
	Let $G$ be a connected graph with $n\geq2$ vertices, $A=A(G)$ and $\rho$ its spectral radius. Then:
	\begin{enumerate}[(i)]
		\item $\rho$ is a simple eigenvalue of $A$, with an eigenvector $v$ such that $v_i>0$ for every $i$;
		\item every eigenvalue $\lambda$ of $A$ satisfies $|\lambda|\leq\rho$, and $-\rho$ is an eigenvalue of $A$ if and only if $G$ is bipartite;
		\item If $G$ is bipartite with classes $V_1\cup V_2$ and $D=\Diag(\varepsilon_1,\dots,\varepsilon_n)$, where $\varepsilon_j=1$ if $j\in V_1$ and $\varepsilon_j=-1$ if $j\in V_2$, then $DAD=-A$. Consequently $-\rho$ is also simple, with eigenvector $Dv$.
	\end{enumerate}
	In particular, $\|A\|=\rho$.
\end{lemma}

\begin{proof}
	The matrix $A$ is nonnegative, and it is irreducible because $G$ is connected \cite[Theorem~6.2.24]{HJ}. Then (i) and the first part of (ii) follow from the Perron--Frobenius theorem for irreducible nonnegative matrices \cite[Theorem~8.4.4]{HJ}, and the characterization of bipartite graphs in (ii) is \cite[Proposition~3.4.1]{BH}. For (iii), $(DAD)_{ij}=\varepsilon_i\varepsilon_j a_{ij}$ and every edge joins different classes; since $D=D^{-1}$, the similarity $DAD=-A$ maps the eigenspace of $\rho$ onto that of $-\rho$.
\end{proof}

\begin{proposition}\label{prop:grafos}
	Let $G$ be a connected graph and $A_0=A(G)$. Then $A_0$ is $\calW$-minimal if and only if $G$ is complete bipartite. In that case, $\text{dist}(A_0,\calW)=\|A_0\|=\sqrt{mk}$ for $G=K_{m,k}$.
\end{proposition}

\begin{proof}
	By the equivalence of items (1) and (5) of Theorem \ref{teo minimality merged}, $A_0$ is $\calW$-minimal if and only if there exists $F\in\partial\|A_0\|\cap\calW^\circ_1$, and such an $F$ has zero diagonal. Let $\rho=\|A_0\|$ and let $v$ be the unit eigenvector of Lemma \ref{lem:PF}\,(i).

	If $G$ is not bipartite, by Lemma \ref{lem:PF}\,(i) and (ii), $\lambda_{\max}=\rho>-\lambda_{\min}$ and $\rho$ is simple. By \eqref{subdif no balanceado}, $\partial\|A_0\|=\{vv^*\}$, whose diagonal entries $v_i^2$ are positive. Hence $A_0$ is not $\calW$-minimal.

	If $G$ is bipartite, by Lemma \ref{lem:PF}\,(iii), $\pm\rho$ are simple eigenvalues with unit eigenvectors $v$ and $w=Dv$, that is, $w_j=\varepsilon_jv_j$. By \eqref{subdif 2} in Lemma \ref{lema of norm subdif}, the Hermitian elements of $\partial\|A_0\|$ are
	\[F=\alpha\, vv^*-(1-\alpha)\, ww^*,\qquad \alpha\in[0,1].\]
	The diagonal of such an $F$ is $\left( (2\alpha-1)v_i^2\right)_i$, so that a zero diagonal forces $\alpha=\frac12$. For this $\alpha$,
	\[F_{ij}=\tfrac12\, v_iv_j\,(1-\varepsilon_i\varepsilon_j)=
	\begin{cases}
		0 & \text{if } i,j \text{ belong to the same class},\\
		v_iv_j\neq0 & \text{if } i,j \text{ belong to different classes}.
	\end{cases}\]
	Therefore $F\in\calW^\circ_1$ if and only if every pair of vertices in different classes is an edge, that is, if and only if $G$ is complete bipartite. For $G=K_{m,k}$ we have $\rho=\sqrt{mk}$.
\end{proof}

If the vertices of $K_{m,k}$ are numbered so that the class of size $m$ comes first, then $A_0=\begin{psmallmatrix}0&J\\J^T&0\end{psmallmatrix}$, where $J$ is the $m\times k$ matrix of ones, and $\calW=M_m^h(\C)\oplus M_k^h(\C)$ is the Hermitian part of a unital C$^*$-subalgebra, so this case is covered by \cite{BV-2026}. Here the witness of the proof is $F=\frac12(vv^*-ww^*)=A_0/\|A_0\|_1$. Note also that, by Lemma \ref{lem:PF}\,(ii), the spectrum of $A(G)$ is balanced exactly when $G$ is bipartite, so bipartite graphs fall into the second condition of \eqref{0 in moments} and non-bipartite graphs into the first one; both conditions are needed.

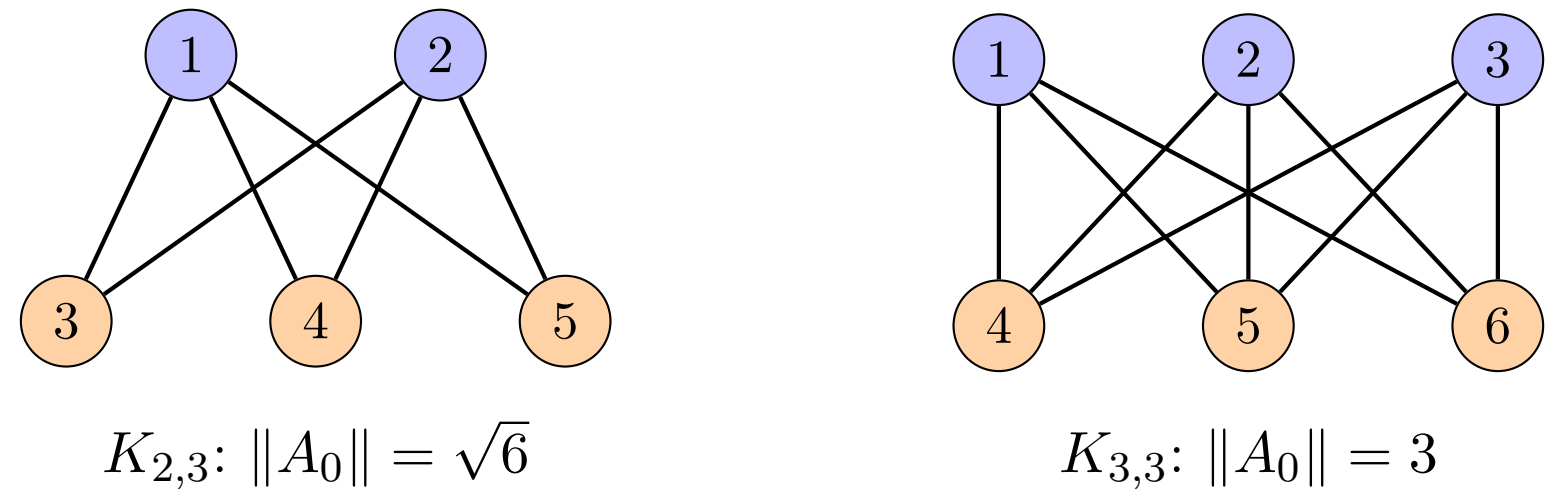
\begin{figure}[ht]
	\centering
	\begin{tikzpicture}[scale=1.1]
		% K_{2,3}
		\foreach \i/\x in {1/0.75,2/2.25} \node[vA] (a\i) at (\x,1.6) {\i};
		\foreach \j/\x in {3/0,4/1.5,5/3} \node[vB] (b\j) at (\x,0) {\j};
		\foreach \i in {1,2} \foreach \j in {3,4,5} \draw[arista] (a\i) -- (b\j);
		\node at (1.5,-0.8) {$K_{2,3}$:\ $\|A_0\|=\sqrt6$};
	\end{tikzpicture}
	\hspace{2cm}
	\begin{tikzpicture}[scale=1.1]
		% K_{3,3}
		\foreach \i/\x in {1/0,2/1.5,3/3} \node[vA] (a\i) at (\x,1.6) {\i};
		\foreach \j/\x in {4/0,5/1.5,6/3} \node[vB] (b\j) at (\x,0) {\j};
		\foreach \i in {1,2,3} \foreach \j in {4,5,6} \draw[arista] (a\i) -- (b\j);
		\node at (1.5,-0.8) {$K_{3,3}$:\ $\|A_0\|=3$};
	\end{tikzpicture}
	\caption{Two complete bipartite graphs, with the classes of the bipartition in different colors. By Proposition~\ref{prop:grafos} their adjacency matrices are $\calW$-minimal: non-edges always join vertices of the same class, which is exactly where the witness $F=\frac12(vv^*-ww^*)$ vanishes.}
	\label{fig:bipartitos}
\end{figure}

\begin{figure}[ht]
	\centering
	\begin{tikzpicture}[scale=1.1]
		\node[vA] (1) at (0,0) {1};
		\node[vB] (2) at (1.5,0) {2};
		\node[vA] (3) at (3,0) {3};
		\node[vB] (4) at (4.5,0) {4};
		\draw[arista] (1) -- (2) -- (3) -- (4);
		\node at (2.25,-1.9) {$A_0=A(P_4)$:\ $\|A_0\|=\varphi$};
	\end{tikzpicture}
	\hspace{1.5cm}
	\begin{tikzpicture}[scale=1.1]
		\node[vA] (1) at (0,0) {1};
		\node[vB] (2) at (1.5,0) {2};
		\node[vA] (3) at (3,0) {3};
		\node[vB] (4) at (4.5,0) {4};
		\draw[arista] (1) -- (2) -- (3) -- (4);
		\draw[noarista] (1) to[bend right=40] node[below,font=\small] {$-1$} (4);
		\node at (2.25,-1.9) {$A_0+W_0$:\ $\|A_0+W_0\|=\sqrt2$};
	\end{tikzpicture}
	\caption{The path $P_4$ and its optimal completion. The pair $\{1,4\}$ joins different classes but is not an edge, and this is exactly where the witness of Proposition~\ref{prop:grafos} fails. The optimal completion has weight $-1$ at this position (dashed).}
	\label{fig:P4}
\end{figure}
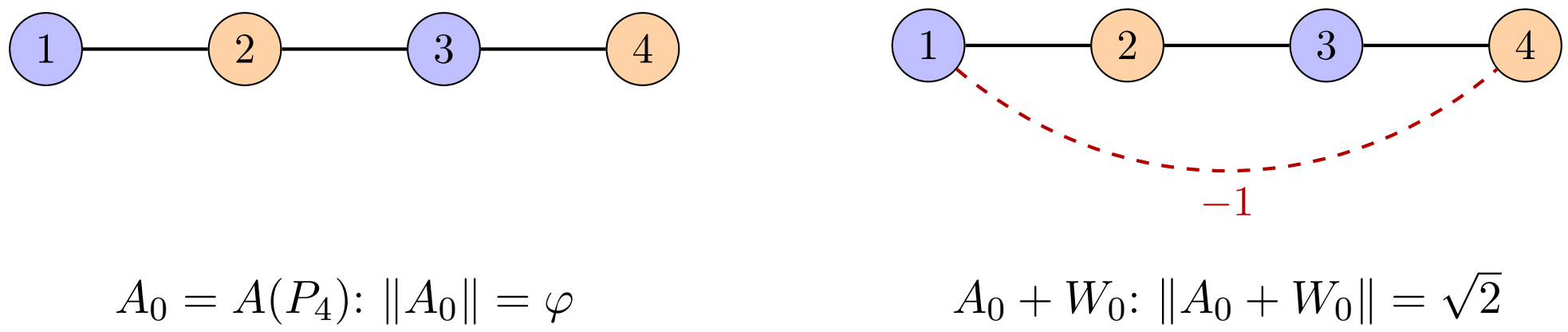

\begin{example}\label{ex:P4}
	Completeness cannot be dropped in Proposition \ref{prop:grafos}. Let $A_0=A(P_4)$, where $P_4$ is the path $1-2-3-4$, which is bipartite but not complete. Then $\|A_0\|=\varphi=\frac{1+\sqrt5}{2}$, while $\text{dist}(A_0,\calW)=\sqrt2$, attained at $W_0=-\sqrt2\,B^{1,4}$, that is, with weight $-1$ at the non-edge $\{1,4\}$ (see Figure \ref{fig:P4}):
	\[A_0+W_0=\begin{pmatrix}
		0 & 1 & 0 & -1\\
		1 & 0 & 1 & 0\\
		0 & 1 & 0 & 1\\
		-1 & 0 & 1 & 0
	\end{pmatrix},
	\qquad (A_0+W_0)^2=2I.\]
	In particular $\sigma(A_0+W_0)=\{\pm\sqrt2\}$ and $S_{\mathrm{MAX}}=\C^4$. A witness is the matrix $F\in\calW^\circ_1$ supported on the edges, with weights $a=c=\frac{1}{4\sqrt2}$ at $\{1,2\}$, $\{3,4\}$ and $b=\frac{1}{2\sqrt2}$ at $\{2,3\}$. Indeed, for such tridiagonal matrices with $a,c\geq0$ one has $\|F\|_1=2\sqrt{(a+c)^2+b^2}=1$, so that for every $W\in\calW$
	\[\sqrt2=2(a+b+c)=\tr(FA_0)=\tr(F(A_0+W))\leq\|A_0+W\|,\]
	with equality at $W_0$. Every choice with $a,c\geq0$ and $a+c=b=\frac{1}{2\sqrt2}$ works, so the dual witness is not unique.

	The minimality depends on $\calW$: if only the diagonal is free, that is, $I_0=\{(i,i)\}_{1\leq i\leq 4}$, then $A(P_4)$ is minimal, since it is tridiagonal with zero diagonal \cite[Lemma~6.3]{B-21}, and $\text{dist}(A(P_4),\Diag(M_4^h(\C)))=\varphi>\sqrt2$. More generally, $A(G)$ is minimal with respect to the real diagonal matrices for every connected bipartite graph $G$, because the matrix $F=\frac12(vv^*-ww^*)$ of the proof of Proposition \ref{prop:grafos} always has zero diagonal; what fails outside complete bipartite graphs is only that $F$ does not vanish at the non-edges joining different classes.
\end{example}

\begin{example}\label{ex:Kn}
	For $G=K_n$ there are no non-edges, so $\calW$ is the space of real diagonal matrices, which is the case studied in \cite{alrv}. For $n\geq3$ the graph is not bipartite, and indeed $\|A(K_n)\|=n-1$ while
	\[\text{dist}\left( A(K_n),\calW\right) =\frac n2,\]
	attained at $J-\frac n2 I$, where $J$ is the $n\times n$ matrix of ones; its spectrum is $\left\{ \frac n2,-\frac n2\right\}$, with $-\frac n2$ of multiplicity $n-1$. With $v=\frac{1}{\sqrt n}\mathbf 1$ and $P_{\mathbf 1^\perp}=I-vv^*$, a witness is
	\[F=\tfrac12\left( vv^*-\tfrac{1}{n-1}P_{\mathbf 1^\perp}\right),\]
	which has zero diagonal, because both terms have constant diagonal equal to $\frac{1}{2n}$, and satisfies $\|F\|_1=1$ and $\tr\left( F(J-\tfrac n2I)\right)=\frac n2$.
\end{example}

\end{document}